\documentclass[11pt, letterpaper]{amsart}
\usepackage{graphicx, amssymb}
\usepackage{amsmath}

\newtheorem{thm}{Theorem}[section]
\newtheorem{cor}[thm]{Corollary}
\newtheorem{lem}[thm]{Lemma}
\newtheorem{prop}[thm]{Proposition}
\theoremstyle{definition}
\newtheorem{defn}[thm]{Definition}

\theoremstyle{remark}
\newtheorem{rem}[thm]{Remark}
\newtheorem{exa}[thm]{Example}

\numberwithin{equation}{section}

\newcommand{\set}[1]{\left\{#1\right\}}
\newcommand{\Real}{\mathbb R}
\newcommand{\Natural}{\mathbb N}

\newcommand{\F}{\mathcal{F}}

\newcommand{\T}{\mathcal{T}}
\newcommand{\prob}{\mathbb{P}}

\newcommand{\C}{\mathcal{C}}
\newcommand{\expec}{\mathbb{E}}

\newcommand{\cS}{\mathcal{S}}
\newcommand{\cM}{\mathcal{M}}

\newcommand{\bL}{\mathbb{L}}
\newcommand{\cL}{\mathcal{L}}
\newcommand{\cadlag}{c\`adl\`ag}
\newcommand{\such}{\, : \, }

\newcommand{\basis}{(\Omega,  \, (\F_t)_{t \in \Real_+}, \, \prob)}
\newcommand{\indic}{\mathbb{I}}
\newcommand{\pare}[1]{\left(#1\right)}
\newcommand{\bra}[1]{\left[#1\right]}
\newcommand{\wt}[1]{\widetilde{#1}}
\newcommand{\vs}{\vspace*{2mm}}

\newcommand{\oY}{\overline{Y}}
\newcommand{\oZ}{\overline{Z}}
\newcommand{\cY}{\mathcal{Y}}
\newcommand{\cZ}{\mathcal{Z}}
\newcommand{\kY}{{}^k \kern-0.15em Y}
\newcommand{\kZ}{{}^k \kern-0.15em Z}
\newcommand{\uX}{\underline{X}}
\newcommand{\bY}{\mathbb{Y}}
\newcommand{\bZ}{\mathbb{Z}}
\newcommand{\og}{\overline{g}}
\newcommand{\ou}{\overline{u}}

\title[BSDE and strict local martingale]{On Backward stochastic differential equations and \\strict local martingales}

\author[]{Hao Xing}
\thanks{This research is supported in part by STICERD at London School of Economics. We are grateful to two anonymous referees and the Associate Editor for their valuable comments, which helped us improve this paper.}
\address{Department of Statistics, London School of Economics and Political Science, 10 Houghton st, London, WC2A 2AE, UK}
\email{h.xing@lse.ac.uk}

\date{}

\begin{document}


\begin{abstract}
We study a backward stochastic differential equation whose terminal condition is an integrable function of a local martingale and generator has bounded growth in $z$. When the local martingale is a strict local martingale, the BSDE admits at least two different solutions. Other than a solution whose first component is of class \emph{D}, there exists another solution whose first component is not of class \emph{D} and strictly dominates the class \emph{D} solution. Both solutions are $\mathbb{L}^p$ integrable for any $0<p<1$. These two different BSDE solutions generate different viscosity solutions to the associated quasi-linear partial differential equation. On the contrary, when a Lyapunov function exists, the local martingale is a martingale and the quasi-linear equation admits a unique viscosity solution of at most linear growth.

\end{abstract}

\date{December 5, 2011}

\maketitle

\setcounter{section}{-1}

\section{Introduction}\label{sec: intro}

Let $B=\set{B_t \such t\geq 0}$ be a standard $d$-dimensional Brownian motion defined on some complete probability space $\basis$. Here $\{\F_t\}_{t\geq 0}$ is the argumented natural filtration of $B$ which satisfies the natural conditions. Fix a real number $T>0$. Consider a continuous adapted process $\{X_t: t\in [0,T]\}$ on $\basis$ with value in $\Real_+^d$ such that each component of $X$ is a nonnegative local martingale. Here $X$ may not necessarily be Markovian. We call $X$ a \emph{martingale}, if all its components are martingales, otherwise $X$ is a \emph{strict local martingales}.

Given a terminal function $g:\Real_+^d \rightarrow \Real$ and a generator $f: [0,T]\times \Real_+^d \times \Real\times \Real^d \rightarrow \Real$, we consider the following backward stochastic differential equation:
\begin{equation}\label{eq: bsde}
 Y_t = g(X_T) + \int_t^T f(s, X_s, Y_s, Z_s) \, ds - \int_t^T Z_s \, dB_s, \quad 0\leq t\leq T. \tag{BSDE}
\end{equation}
We look for progressively measurable processes $(Y,Z) = \set{(Y_t, Z_t) \such t\in [0,T]}$ such that they satisfy the previous equation $\prob$-a.s. and every term in the equation is well defined. Such equation, in the nonlinear case, is a special type of backward stochastic differential equations (BSDE) introduced in \cite{Pardoux-Peng}. Since then, BSDEs have been studied with great interest.

Let us briefly review existence and uniqueness results for BSDE solutions with different integrability properties. When $g(X_T)$ and $\set{f(t, X_t, 0, 0) \such t\in [0,T]}$, which are called \emph{parameters}, are square integrable, Pardoux and Peng proved in \cite{Pardoux-Peng} the existence and uniqueness for the square integrable ($\bL^2$-) solution of BSDEs with Lipschitz continuous generators. When parameters are $\bL^p$ with $p\in (1,2)$, the existence of $\bL^p$-solutions was established by El Karoui et al. in \cite{El-Karoui-et-al}, and later extended by Briand et al. in \cite{Briand-et-al}, where a uniqueness result was also obtained. For only $\bL^1$-integrable parameters, Peng  studied a BSDE  in \cite{Peng-97} whose generator is a sum of two functions in $y$ and $z$ respectively. This was extended to BSDEs whose generator has strictly sublinear growth in $z$ by Briand et al. in \cite{Briand-et-al}. In this paper, existence and uniqueness of solutions have been established in \emph{class D}, i.e., the class of processes $Y$ such that $\set{Y_\tau \such \tau \text{ is } \F- \text{stopping time with value in } [0,T]}$ is uniformly integrable. However, all the above results do not cover the following example, which motivates this study.

Consider the following stochastic differential equation (SDE):
\begin{equation}\label{eq: inv-bessel}
 dX_t = - X_t^2 \, dW_t, \quad X_0= x >0,
\end{equation}
where $W$ is a standard $1$-dimensional Brownian motion. This SDE admits a unique nonnegative strong solution $\{X_t \such t\geq 0\}$, which is the so called reciprocal $3$-dimensional Bessel process. It is well known that $X$ is a strict local martingale and $\expec[X_T^2]<\infty$ (see (2.13) in \cite{Revuz-Yor} pp. 194). Let us consider the following BSDE with zero generator:
\begin{equation}\label{eq: BSDE inv-bessel}
 Y_t = X_T - \int_t^T Z_s \, dW_s, \quad 0\leq t\leq T.
\end{equation}
It follows from the martingale representation theorem that $\oY_\cdot = \expec[X_T \,|\, \F_\cdot]$ and its associated integrand $\oZ$ solve the previous BSDE. Moreover the Burkholder-Davis-Gundy inequality (see e.g. Theorem 42.1 in \cite{Rogers-Williams} Chap. IV) implies that both $\expec[\sup_{0\leq t\leq T} \oY^2_t]$ and $\expec[\int_0^T \oZ_s^2\, ds]$ are finite. Therefore $(\oY, \oZ)$ is an $\bL^2$-solution.

However, there is another obvious solution to \eqref{eq: BSDE inv-bessel}. That is $(Y, Z) =  (X, -X^2)$. To the best of our knowledge, this solution has not been studied before. This solution solves \eqref{eq: BSDE inv-bessel}, but it does not satisfy integrability properties reviewed earlier. To begin with, $\expec[\int_0^T Z_s^2 \,ds] = \expec[\int_0^T (X_s^2)^2\, ds] =\infty$. If $X^2$ was square-integrable, $\int_0^\cdot X_s^2 \,dW_s$ would be a martingale. This implies $X_0=\expec[X_T]$ which contradicts with the strict local martingale property of $X$. Additionally, $Y=X$ is clearly not of class $D$, otherwise $X$ would be a martingale again. Moreover $\expec[\sup_{0\leq t\leq T} Y_t] = \expec[\sup_{0\leq t\leq T} X_t] =\infty$, which implies $\expec[(\int_0^T Z_s^2 \,ds)^{1/2}] = \expec[(\int_0^T (X_s^2)^2 \,ds)^{1/2}] =\infty$ from the Burkholder-Davis-Gundy inequality.

Nevertheless Lemma \ref{lem: local mart Lp} below shows that $\expec[\sup_{0\leq t\leq T} Y_t^p] <\infty$ for any $p\in (0,1)$.  Hence $\expec[(\int_0^T Z_s^2 \,ds)^{p/2}] <\infty$ follows from the Burkholder-Davis-Gundy inequality. Therefore $(Y, Z)$ is one $\bL^p$ $(p\in (0,1))$ solution to \eqref{eq: BSDE inv-bessel}. On the other hand, Jensen's inequality entails that $(\oY, \oZ)$ is also an $\bL^p$ solution. Therefore there are at least two solutions to \eqref{eq: BSDE inv-bessel} inside the same class of processes.

The previous example is closely related to the notion of $g$-martingales introduced in \cite{Peng-97}. The BSDE solutions can be considered as nonlinear martingales because a solution to BSDE with zero generator is given by conditional expectation of the terminal condition. In classical theory, martingales are local martingales. Therefore to have a nonlinear theory which contains the classical theory, it is necessary to extend the notion of local martingales into the framework of BSDEs. In this paper, we regard solutions to \eqref{eq: bsde} as \emph{$g$-local martingales}. When $X$ is a classical strict local martingale, other than the class \emph{D} solution obtained in \cite{Briand-et-al}, there exists another solution which is not of class \emph{D}. We regard it as a \emph{$g$-strict local martingale}. Example in \eqref{eq: BSDE inv-bessel} is a special example of \eqref{eq: bsde}.

Another motivation of this paper is to study the connection between \eqref{eq: bsde} and its associated quasi-linear partial differential equation (PDE). When $X$ is a diffusion whose dynamics is $dX_t = \sigma(X_t) \,dB_t$, the quasi-linear PDE associated to \eqref{eq: bsde} reads
\begin{equation}\label{eq: pde}
\tag{PDE}
\begin{array}{ll}
 - \partial_t u - \frac12 Tr\pare{\sigma \sigma' \nabla^2 u} - f(t, x, u, \nabla u \,\sigma) =0, &\quad (t,x) \in [0,T) \times (0,\infty)^d,\\
 u(T, x) = g(x), & \quad x\in (0, \infty)^d.
\end{array}
\end{equation}
Since the dawn of the BSDE theory, close connections between BSDEs and quasi-linear PDEs have been established (see e.g. \cite{Pardoux-Peng-PDE} and \cite{Barles-et-al}). These results may be seen as generalizations of the celebrated Feynman-Kac formula. Since \eqref{eq: bsde} may have multiple solutions, it is natural to expect multiple solutions to \eqref{eq: pde}. Actually, when $f$ vanishes, $g$ has linear growth, and $X$ is a strict local martingale, multiple solutions to \eqref{eq: pde} (now a linear equation) has been observed in \cite{Heston-Loewenstein-Willard}. See \cite{Ekstrom-Tysk}, \cite{bx09} and \cite{bkx10} for recent developments. In these studies, X being a martingale has been shown to be the necessary and sufficient condition for the uniqueness of classical solutions, in the class of at most linear growth functions, to valuation equations associated with local/stochastic volatility models. However existing results treat PDEs with $1$ or $2$ spatial dimension and employ the notion of classical solutions. When the equation is nonlinear, classical solutions are in general not expected. It is then natural to work in the framework of viscosity solutions. However when $X$ is a strict local martingale, its volatility coefficient $\sigma$ fails to be Lipschitz on the entire state space. Therefore classical techniques in viscosity solutions need to be extended to treat local Lipschitz coefficients. See \cite{Amadori} and \cite{Constantini-et-al} for recent developments in this direction.
\vspace{0.2cm}

{\bf Our work:} Assume that $g$ is nonnegative and has at most linear growth, $f$ satisfies a monotonicity condition in $y$ and has bounded growth $z$. When $X$ is a strict local martingale, \eqref{eq: bsde} admits at least two solutions. The first component of one solution is of class \emph{D}. Theorem \ref{thm: bsde existence} shows that there exists another solution whose first component is not of class \emph{D} and is strictly larger than the class \emph{D} solution. These two BSDE solutions induce different viscosity solutions to \eqref{eq: pde}. See Theorem \ref{thm: pde existence}. On the other hand, when a Lyapunov function exists, $X$ is a martingale, moreover Theorem \ref{thm: pde uniqueness} shows that \eqref{eq: pde} admits a unique viscosity solution in the class of functions with at most linear growth. Contrast to the existing results on the uniqueness of viscosity solutions for PDEs with global Lipschitz coefficients, the volatility coefficient of $X$ is assumed to be only locally Lipschitz continuous.

\vspace{0.2cm}

Multiple solutions of BSDEs have been observed by Bao et al. in \cite{Bao-Delbaen-Hu}. Contrast to their source of multiplicity, which is the multiple choices of boundary conditions for the associated PDE, our multiple solutions are induced by the linear growth terminal condition and the strict local martingale property of $X$. When $X$ does not explode to the boundary of its state space, no boundary condition is needed for \eqref{eq: pde}, multiple solutions still exist (see Theorem \ref{thm: pde existence}).

Even though the generator $f$ is assumed to have bounded growth in $z$, \eqref{eq: bsde} is related to some special quadratic BSDEs, whose generator has quadratic growth in $z$, via the exponential transform. As a result, explicit multiple solutions to these quadratic BSDEs are constructed in Example \ref{exa: quad growth in z}. We refer readers to \cite{Kobylanski}, \cite{Briand-Hu}, and \cite{Barrieu-El-Karoui} for existence results of solutions to quadratic BSDEs and \cite{Delbaen-Hu-Richou} for uniqueness results.

The rest of the paper is organized as follows. After notation and definitions are introduced, we present our main results in Section \ref{sec: main results}. Several examples are given in this section to illustrate our results.  Multiple BSDE solutions are constructed in Section \ref{sec: bsde solution}. Existence and uniqueness of viscosity solutions are proved in Section \ref{sec: pde solution}.

\section{Main results}\label{sec: main results}
\subsection{Notation and definitions}\label{subsec: notation}

Throughout this paper, we fix the probability measure $\prob$. Every relationship between random variables is understood in $\prob$-almost sure sense.

For any $p>0$, $\cS^p$ denotes the class of real valued, adapted and \cadlag\, process $\set{Y_t; t\in [0,T]}$ such that
\[
 \|Y\|_{\cS^p} := \expec[\sup_{t\in[0,T]} |Y_t|^p]^{1 \wedge 1/p} <+\infty.
\]
If $p\geq 1$, $\|\cdot\|_{\cS^p}$ is a norm on $\cS^p$ and if $p\in (0,1)$, $(Y, Y') \mapsto \|Y-Y'\|_{\cS^p}$ denotes a distance on $\cS^p$. Under this metric, $\cS^p$ is complete. We denote $\cS^\infty$ the set of adapted bounded processes. Denote by $\T_{[0,T]}$ the set of all $\F$-stopping time $\tau$ such that $0\leq \tau\leq T$. We call $Y$ belongs to the \emph{class D} if the family $\set{Y_\tau; \tau\in \T_{[0,T]}}$ is uniformly integrable.
Let $\cM^p$ denote the equivalent class of predictable processes $\set{Z_t; t\in [0,T]}$ with values in $\Real^d$ such that
\[
 \|Z\|_{\cM^p}:= \expec\bra{\pare{\int_0^T |Z_s|^2 \,ds}^{p/2}}^{1\wedge 1/p}<+\infty.
\]
For $p\geq 1$, $\cM^p$ is a Banach space with this norm, and for $p\in(0,1)$, $\cM^p$ is a complete metric space with the resulting distance.

The Euclidean norm is denoted as $|\cdot|$ regardless of dimension. Denote $\mathbb{B}_r := \{x\in \Real_+^d: |x|<r\}$, $\mathbb{B}^+_r := \{x\in (0,\infty)^d: |x|<r\}$, and $\mathbb{S}^+_r := \{x\in (0,\infty)^d: |x|=r\}$ for some $r>0$.
For $x\in \Real^d$, $x^i$ is its $i-$component and $\underline{x} := \sum_{i=1}^d x^i$. For the process $X$, we denote \[\uX = \sum_{i=1}^d X^i.\]

Let us recall what we mean by a solution to \eqref{eq: bsde}.
\begin{defn}\label{def: solution}
 A solution to \eqref{eq: bsde} is a pair $(Y,Z) = \{(Y_t, Z_t) : t\in [0,T]\}$ of progressively measurable processes with values in $\Real \times \Real^d$ such that $\prob$-a.s. $t\mapsto Y_t$ is continuous, $\int_0^T |Z_t|^2 dt<\infty$, $\int_0^T |f(t, X_t, Y_t, Z_t)| \,dt <\infty$, and \eqref{eq: bsde} is satisfied.
\end{defn}

\subsection{Existence of BSDE solutions}\label{subsec: existence BSDE soln}

As in the Introduction, each component of $X$ is a nonnegative local martingale. Hence both $X^i$, $1\leq i\leq d$, and $\uX$ are supermartingales. The terminal function $g$ is continuous, nonnegative, and
\begin{equation}\label{H: g linear growth}\tag{H1}
 K:= \sup\left\{ \frac{g(x)}{1+\underline{x}} \such x\in \Real_+^d\right\} <\infty.
\end{equation}
Hence $0\leq g(x) \leq K(1+\underline{x})$ for any $x\in \Real_+^d$.
Combined with the supermartingale property of $\uX$, \eqref{H: g linear growth} implies $g(X_T) \in \bL^1$. Since we focus on only integrable terminal conditions, we do not a priori assume $g(X_T) \in \bL^p$ for some $p>1$. If the parameters are $\bL^p$-integrable for some $p>1$, existence and uniqueness of solutions in $(\cS^p, \cM^p)$ have been established in \cite{Briand-et-al}.

For the generator, we assume that $f$ is jointly continuous in all its variables. Moreover, there exists a function $H: [0,T] \times \Real_+ \rightarrow \Real_+$ such that
\begin{align}
 & H \text{ is locally bounded on } [0,T]\times \Real_+, \label{H: H local bound}\tag{H2.i}\\
 & r\mapsto H(t, r) \text{ is nondecreasing and concave}. \label{H: H nondec concave}\tag{H2.ii}
\end{align}
There exist constants $\nu $ and $\mu$ such that, for each $(t, x, y, y', z, z') \in [0,T] \times \Real_+^d \times \Real\times \Real \times \Real^d \times \Real^d$,
\begin{align}
 & |f(t,x,y,z) - f(t,x,y,z')| \leq \nu |z-z'|, \label{H: f lip z}\tag{H3.i}\\
 & (y-y')(f(t,x,y,z) - f(t,x,y',z)) \leq \mu(y-y')^2, \label{H: f mono y}\tag{H3.ii}\\
 & f(t, x, y, z) \geq 0, \label{H: f nonneg}\tag{H3.iii}\\
 & f(t, x, 0, z) \leq H(t, \underline{x}). \label{H: f growth z}\tag{H3.iv}
\end{align}

\begin{rem}
When $g$ and $f$ only depend on some components of $X$, sums on $X_i$ should be taken only on these components. All results in this paper still hold. For simplicity of presentation, we assume that both $g$ and $f$ depend nontrivially on all components of $X$.
\end{rem}

\begin{rem}
 Since both $g$ and $f$ are nonnegative, one can expect that we are interested to find solutions with nonnegative first component. Assumptions \eqref{H: f mono y} and \eqref{H: f growth z} combined yields that $f(t,x,y,z) \leq \mu y + H(t, \underline{x})$ for $(t,x,y,z) \in [0,T] \times \Real_+^d \times \Real_+ \times \Real^d$. Hence $f$ has bounded growth in $z$. This assumption, together with the assumptions on $H$, will facilitate the construction of \eqref{eq: bsde} solutions and imply that their first component is inside the following class.
\end{rem}

Let us define a class of continuous adapted processes:
\[
\C:= \set{Y \such 0\leq Y_t \leq C \pare{K \pare{1+ \uX_t} + \expec\bra{\left.\int_t^T H(s, \uX_s)\, ds \, \right| \,\F_t}} \text{ for any } t\in [0, T]},
\]
where $C = e^{(\mu \vee 0) T}$.
For a solution $(Y, Z)$ to \eqref{eq: bsde} such that $Y\in \C$, Proposition \ref{prop: integrability} below shows that $(Y, Z)\in (\cS^p, \cM^p)$ for any $p\in (0,1)$. We are now ready to present the first main result.

\begin{thm}\label{thm: bsde existence}
 Let (H1) - (H3) hold.
 \begin{enumerate}
  \item[(i)] There exists a solution $(\oY, \oZ)$ to \eqref{eq: bsde} such that $\oY \in \C$ and $\oY$ is of class \emph{D}.
  \item[(ii)] For any other solution $(\widetilde{Y}, \widetilde{Z})$ to \eqref{eq: bsde} such that $\widetilde{Y} \in \C$, $\widetilde{Y}_t \geq \oY_t$ for any $t\in [0,T]$.
 \end{enumerate}
 Define $\overline{g}(x):= K(1+\underline{x}) - g(x)$. Assume that $\overline{g}$ satisfies the following assumptions:
 \begin{align}
  & \overline{g}(X_\cdot) \text{ is a supermartingale on } [0,T], \label{H: og supermart}\tag{H4.i}\\
  & \text{there exists a nondecreasing univariate continuous function } \overline{G} : \Real_+ \rightarrow \Real_+ \text{ such that} \notag\\
  & \hspace{3cm} \overline{g}(x) \leq \overline{G}(\underline{x}) \quad \text{ and } \quad \lim_{r\rightarrow \infty} \overline{G}(r)/r =0. \label{H: og strict sublinear} \tag{H4.ii}
 \end{align}
 \begin{enumerate}
  \item[(iii)] Then when $X$ is a strict local martingale on $[0,T]$, there exists another solution $(Y, Z)$ such that $Y\in \C$ and $Y_\cdot \geq K\pare{\uX_\cdot - \expec[\uX_T \,|\, \F_\cdot]} + \expec[g(X_T) \,|\, \F_\cdot]$, but $Y$ is not of class \emph{D}, moreover $Y_0 > \oY_0$.
 \end{enumerate}
\end{thm}

\begin{rem}
 The existence of different solutions to the same BSDE implies that the comparison result for BSDE solutions fails in class $\C$. To restore the comparison in $\C$, one can assume
 \[
  K\pare{1+\uX_\cdot} + \expec\bra{\left. \int_\cdot^T H(s, \uX_s) \,ds \right| \,\F_\cdot} \in \bigcup_{p>1} \cS_p.
 \]
 Indeed, this condition yields $Y\in \bigcup_{p>1} \cS_p$ for any $Y\in \C$. Then the comparison result for solutions in class $\C$ follows from Proposition 5 in \cite{Briand-Hu}. It should be pointed out that this condition already excludes strict local martingales $X$.
\end{rem}

\begin{rem}
 The solution whose first component is of class \emph{D} is unique, if the following additional assumption on $f$ is satisfied: there exist two constants $\gamma\geq 0$ and $\beta\in (0,1)$ such that
 \[
  |f(t,x,y,z) - f(t,x,y,0)| \leq \gamma (\mu y + H(t, \underline{x}))^\beta, \quad \text{ for all } (t,x,y,z) \in [0,T] \times \Real_+^d \times \Real_+ \times \Real^d.
 \]
 This follows from Theorem 6.2 in \cite{Briand-et-al}. Note that the above assumption is trivially satisfied if $f$ does not depend on $z$.
\end{rem}

It has been observed in \cite{Ekstrom-Tysk} that linear \eqref{eq: pde} admits an uncountable family of different solutions when $X$ is a strict local martingale. This translates to an uncountable family of different solutions to the associated BSDE which has zero generator. This phenomenon can be extended to BSDEs with nonzero generators as follows.
\begin{cor}\label{cor: inf soln}
 Let (H1) - (H4) hold. Assume that $f$ is Lipschitz in $y$ and does not depend on $z$. Then \eqref{eq: bsde} admits a family of solutions $(Y^\alpha, Z^\alpha)_{\alpha\in [0,1]}$ with $(Y^0, Z^0) = (\overline{Y}, \overline{Z})$ and $(Y^1, Z^1) = (Y, Z)$, moreover $\{Y^\alpha\}_{\alpha \in [0,1]}$ is nondecreasing sequence in $\mathcal{C}$ such that $\{Y^\alpha_0\}_{\alpha\in [0,1]}$ is strictly increasing.
\end{cor}

Let us now illustrate Theorem~\ref{thm: bsde existence} in the following three examples. The first example gives a class of terminal conditions which satisfy (H4). This class contains call option payoffs in financial applications.
In the second example, solutions $Y$ and $\oY$ are constructed explicitly when the generator vanishes. The third example presents multiple solutions to a quadratic BSDE. BSDEs in last two examples actually admit uncountable families of different solutions because their generators satisfy assumptions in the previous corollary.

\begin{exa}[Assumption (H4)]\label{exa: assumption H4}
 Assumptions \eqref{H: og supermart} and \eqref{H: og strict sublinear} hold when $g(x)= G(\underline{x})$ for a convex univariate function $G: \Real_+ \rightarrow \Real_+$ such that $\lim_{r\rightarrow \infty} G(r)/r =K$.

 Indeed, $r \mapsto \overline{g}(r) := K(1+ r) - G(r)$ is a nonnegative nondecreasing concave function, moreover $\lim_{r\rightarrow \infty} \overline{g}(r)/r =0$. Let $\tau_n := \inf\set{t\geq 0 \such X_t \notin \overline{\mathbb{B}_n}} \wedge T$ for $n\geq 0$. This sequence of stopping times localizes each component of $X$ and also $\uX$. Moreover $\lim_{n\rightarrow \infty} \tau_n = T$. It then follows from Fatou's lemma and the concavity of $\overline{g}$ that
 \[
 \begin{split}
  \expec\bra{\overline{g}(\uX_t) \,|\, \F_u} &= \expec\bra{\lim_{n \rightarrow \infty} \overline{g}(\uX_{t\wedge \tau_n}) \,|\, \F_u} \leq \liminf_{n\rightarrow \infty} \expec\bra{\overline{g}(\uX_{t\wedge \tau_n}) \,|\, \F_u} \leq \liminf_{n\rightarrow \infty} \overline{g}\pare{\expec[\uX_{t\wedge \tau_n} \,|\, \F_u]} \\
  &= \liminf_{n\rightarrow \infty} \overline{g}(\uX_{u\wedge \tau_n}) = \overline{g}(\uX_u), \quad \text{ for } 0\leq u\leq t\leq T.
 \end{split}
 \]
 Hence (H4) is satisfied in this case.
\end{exa}

\begin{exa}[Zero generator]\label{exa: zero generator}
Let (H1) and (H4) hold. When the generator $f$ vanishes,
 \[
  \oY_\cdot = \expec\bra{g(X_T) \,|\, \F_\cdot} \quad \text{ and } \quad Y_\cdot = K\pare{\uX_. - \expec[\uX_T \,|\, \F_.]} + \expec[g(X_T) \,|\, \F_\cdot].
 \]
 \eqref{H: g linear growth} yields that $Y_\cdot \leq K(\uX_\cdot - \expec[\uX_T \,|\, \F_\cdot]) + K(1+\expec[\uX_T \,|\, \F_\cdot]) = K(1+\uX_\cdot)$. Therefore both $Y$ and $\oY$ are in $\C$. When $X$ is a strict local martingale on $[0, T]$, it is clear that $\oY_0 = \expec[g(X_T)] < K(\uX_0 - \expec[\uX_T]) + \expec[g(X_T)] = Y_0$. Moreover $\oY$ is of class \emph{D}, but $Y$ is not.
\end{exa}

\begin{exa}[A BSDE with quadratic growth in $z$]\label{exa: quad growth in z}
 Consider the following BSDE:
 \begin{equation}\label{eq: bsde quad in z}
  P_t = \log \uX_T + \int_t^T \pare{\alpha_s + \frac12 |Q_s|^2}\, ds -\int_t^T Q_s \, dB_s,
 \end{equation}
 where $\alpha$ is a nonnegative bounded process. Define $(Y, Z) = (e^P, e^P Q)$. It satisfies
 \begin{equation}\label{eq: transformed bsde}
  Y_t = \uX_T + \int_t^T \alpha_s Y_s \,ds -\int_t^T Z_s \,dB_s.
 \end{equation}
 The previous BSDE satisfies (H1)-(H3).
 When $X$ is a strict local martingale, \eqref{eq: transformed bsde} admits two different solutions, so is \eqref{eq: bsde quad in z}.

 In \cite{Delbaen-Hu-Richou}, the uniqueness of solutions to BSDEs with quadratic growth in $z$ is proved among solutions whose first component $Y$ satisfies
 \[
  \expec\bra{e^{\gamma \sup_{0\leq t\leq T} P_t^+} + e^{\epsilon \sup_{0\leq t\leq T} P_t^-}} <\infty, \quad \text{ for some } \gamma >1 \text{ and } \epsilon >0,
 \]
 where $P^+$ and $P^-$ are positive and negative parts of $P$.
 In this example, the additional solution $(P, Q)$, associated to $(Y, Z)$ in Theorem \ref{thm: bsde existence}, is outside the previous class. Indeed, it follows from Theorem~\ref{thm: bsde existence} (iii) that $e^P \geq \uX$. Then $P^+ \geq \log \max\{\uX, 1\}$, hence
 \[
  \expec\bra{e^{\gamma \sup_{0\leq t\leq T} P_t^+}} \geq \expec[\sup_{0\leq t\leq T} \max\{\uX^\gamma_t, 1\}] \geq \expec [\sup_{0\leq t\leq T} \uX_t^\gamma],
 \]
 where the right-hand-side is infinity for any $\gamma>1$ when $X$ is a strict local martingale.
\end{exa}

\subsection{Existence and Uniqueness of viscosity solutions to a quasi-linear PDE}\label{subsec: pde}

Let us now specify a Markovian dynamics of $X$ and study the quasi-linear PDE associated to \eqref{eq: bsde}. Assume that $\sigma : (0,\infty)^d \rightarrow \Real^{d\times d}$ is locally Lipschitz in $(0,\infty)^d$, i.e., for any compact domain $D\subset (0, \infty)^d$, there exists a constant $L_D$ such that $|\sigma(x) - \sigma(y)| \leq L_D|x-y|$ for any $x,y\in D$. We consider the following SDE:
\begin{equation}\label{eq: sde X}
 dX^{i}_s = \sum_{j=1}^d \sigma_{ij}(X_s) \, dB^j_s, \quad X_0 = x\in (0,\infty)^d, \quad i=1, \cdots, d.
\end{equation}
It is well know that \eqref{eq: sde X} admits a unique strong solution $X^x$ up to its explosion time $\zeta$. Let $\{D_n\}_{n\geq 0}$ be a sequence of bounded open domains such that $\overline{D_n} \subset D_{n+1}$ for $n\geq 0$, and $\bigcup_{n\geq 0} D_n = (0,\infty)^d$. Define $\sigma^x_n := \inf\{t\geq 0 \such X^x_t \notin \overline{D_n}\}$.
Then $\zeta = \lim_{n\rightarrow \infty} \sigma^x_n$. We assume that
\begin{equation}\label{H: X non explode}\tag{H5}
 \prob(\zeta = \infty) =1.
\end{equation}
The assumption above implies that \eqref{eq: sde X} admits a unique $(0,\infty)^d$ valued strong solution $\{X^x_t \such t\geq 0\}$. We denote by $\cL:= \frac12 Tr(\sigma \sigma' \nabla^2)$ its infinitesimal generator.

Since components of $X^x$ are continuous supermartingales, $\zeta=\sigma^x_\infty$. Here $\sigma^x_\infty := \inf\{t\geq 0 \such X^x_t \in \mathcal{O}\}$ where $\mathcal{O}:= \{x\in \Real_+^d \such x_i =0 \text{ for some } i\in \{1, \cdots, d\} \}$ is the face of the first orthant. Therefore \eqref{H: X non explode} is equivalent to $\sigma^x_\infty=\infty$, hence $X^x$ never reaches the boundaries of its state space in finite time. As a result no boundary condition is needed for \eqref{eq: pde}. Still Theorem \ref{thm: pde existence} below shows that \eqref{eq: pde} admits multiple solutions. We refer readers to \cite{bkx10} for a detailed discussion on boundary conditions in stochastic volatility models where the volatility process can reach the boundary of its state space.

\begin{rem}\label{rem: not reach O}
 There are several ways to check whether $\sigma^x_\infty$ is almost surely infinite.

 First, if there exists a Lyapunov function $\Psi$ on $(0,\infty)^d$ such that $\lim_{x \rightarrow \overline{x}} \Psi(x) =\infty$ for any $\overline{x} \in \mathcal{O}$  and a positive constant $\lambda$ such that $\cL \Psi(x) \leq \lambda \Psi(x)$ for any $x\in (0,\infty)^d$, then $\sigma^x_\infty=\infty$ (see Theorem 6.7.1 in \cite{Pinsky}). Second, if $\sigma_{ij}(\cdot)$ is continuously differentiable in $[0,\infty)^d$, the matrix $\sigma \sigma'$ degenerates on $\mathcal{O}$, and the so called \emph{Fichera} drifts $\mathfrak{f}_i(x) = -\frac12 \sum_{j=1}^d \partial_{x_j} (\sigma \sigma')_{ij}(x)$ are nonnegative on $\{x_i=0\}$, for each $i=1, \cdots, d$, then $\sigma^x_\infty=\infty$ (see Theorem 9.4.1 and Corollary 9.4.2 in \cite{friedman-stochastic}). Third, if $\sigma_{ij}(x) = x^i \mathfrak{s}_{ij}(x)$ for some matrix $\mathfrak{s}$, then $Z^i = \log X^i$ satisfies
 $
  dZ^i_t = -\frac12 \sum_{j=1}^d \mathfrak{s}_{ij}^2 (e^{Z_t}) \,dt + \sum_{j=1}^d \mathfrak{s}_{ij} (e^{Z_t}) \, dB^j_t,
 $
 where $e^Z = (e^{Z^1}, \cdots, e^{Z^d})$. Since $|X|$ does not explode to infinity in finite time, then $X$ does not hit $\mathcal{O}$ in finite time if and only if $|Z|$ does not explode to infinity in finite time. Then any sufficient condition which ensures the nonexplosion of $Z$ implies $\sigma^x_\infty=\infty$. For example, Khasminskii provided such a sufficient condition (see e.g. Theorem 52.1 in \cite{Rogers-Williams} Chap. V).
 In 1 dimension, $\sigma^x_\infty=\infty$ can be identified via Feller's test.
\end{rem}

Since no growth assumption is imposed on $\sigma$, $X$ can be strict local martingale. The following are some examples.
\begin{exa}~
\begin{enumerate}
  \item[i)] If there exists some component of $X$, say $X^i$, such that $dX^i_t = \sigma_{ii}(X_t^i) dB^i_t$ where $\sigma\neq 0$ on $(0,\infty)$ and $\int_1^\infty r/\sigma^2_{ii}(r) \,dr <\infty$, then $X^i$, hence $X$, is a strict local martingale (see \cite{Delbaen-Shirakawa}).

  \item[ii)] Suppose that $X$ has the following dynamics
  \[
   dX^1_t = X^1_t \,X^2_t \,dB^1_t,  \qquad dX^2_t = X^2_t (\rho \,dB^1_t + \sqrt{1-\rho^2} \,dB^2_t).
  \]
  Then $X^1$, hence $X$, is a strict local martingale if and only if $\rho>0$ (see \cite{Sin}).

  \item[iii)] A large class of multi-variate local martingales $X$ is provided in \emph{stochastic portfolio theory}, where $X$ models the deflated stock prices. When the market price of risk exists and  there is arbitrage relative to the market, $X$ is a strict local martingale (see Section 6 in \cite{Fernholz-Karatzas-survey} for more details).
\end{enumerate}
\end{exa}

After the dynamics of $X$ is introduced, let us consider \eqref{eq: pde} associated to \eqref{eq: bsde}. The following definition of viscosity solutions follows from \cite{Barles}. For a function $u$ defined on $[0,T]\times (0,\infty)^d$, we denote by $u^*$ (resp. $u_*$) the \emph{upper}- (resp. \emph{lower}-) \emph{semicontinuous envelope} of $u$: for all $(t,x) \in [0,T] \times (0,\infty)^d$,
\[
 u^*(t,x) :=\limsup_{[0,T] \times (0,\infty)^d \ni (t',x') \rightarrow (t,x)} u(t', x') \quad \text{ and } \quad u_*(t,x) :=\liminf_{[0,T] \times (0,\infty)^d \ni (t',x') \rightarrow (t,x)} u(t', x').
\]

\vspace{0.1cm}
\begin{defn}[Viscosity solution]
~
 \begin{itemize}
  \item $u$ is called a \emph{viscosity subsolution} of \eqref{eq: pde} if $u^*<\infty$ on $[0,T]\times (0,\infty)^d$ and if for any $\phi\in C^{1,2}((0,T) \times (0,\infty)^d)$ and $(t,x) \in (0,T)\times (0,\infty)^d$, such that
      $0 = (u^*-\phi)(t, x) \geq (u^* -\phi)(\wt{t}, \wt{x})$ for any $(\wt{t}, \wt{x}) \in(0,T)\times (0,\infty)^d$,
      \[
       -\partial_t \phi - \cL \phi - f(t,x, u^*, \nabla \phi \,\sigma) \leq 0.
      \]
      Moreover, $u^*(T,x) \leq g(x)$ for $x\in (0,\infty)^d$.
  \item The \emph{viscosity supersolution} is defined similarly using $u_*$.
  \item $u$ is called a \emph{viscosity solution} of \eqref{eq: pde} if it is both viscosity sub- and supersolution.
 \end{itemize}
\end{defn}

In what follows viscosity solutions of \eqref{eq: pde} are constructed via solutions of \eqref{eq: bsde}. Since there are multiple solutions of \eqref{eq: bsde}, \eqref{eq: pde} also admits multiple viscosity solutions. For a fixed $t\in [0,T]$, under assumptions of Theorem \ref{thm: bsde existence}, the BSDE
\begin{equation}\label{eq: bsde cY}
 Y_s = g(X^x_{T-t}) + \int_s^{T-t} f(u+t, X^x_u, Y_u, Z_u) \, du -\int_s^{T-t} Z_u \,dB_s,\quad s\in [0, T-t],
\end{equation}
admits two solutions which are denoted by $(\cY^{t,x}, \cZ^{t,x})$ and $(\overline{\cY}^{t,x}, \overline{\cZ}^{t,x})$. Define $B^t_s := B_{(s-t)^+}$ and $X^{t,x}_s := X^x_{(s-t)^+}$ for $s\in [0,T]$. Then $B^t$ is a Brownian motion in its own natural filtration and $X^{t,x}$ is the unique strong solution to
\begin{equation}\label{eq: sde X tx}
 dX^{t,x}_s = \sigma(X^{t,x}_s) \, dB^t_s, \quad X^{t,x}_t = x.
\end{equation}
Define $(Y^{t,x}_s, Z^{t,x}_s):= (\cY^{t,x}_{(s-t)^+}, \indic_{\{s\geq t\}} \cZ^{t,x}_{s-t})$,  and $(\oY^{t,x}_s, \oZ^{t,x}_s)$ similarly, for $s\in [0,T]$. Now two deterministic functions can be defined:
\begin{equation}\label{eq: u and ou}
 u(t,x) := Y^{t,x}_t = \cY^{t,x}_0 \quad \text{ and } \quad \overline{u}(t,x) := \oY^{t,x}_t = \overline{\cY}^{t,x}_0, \quad \text{ for } (t,x)\in [0,T] \times (0,\infty)^d.
\end{equation}
It is immediate from Theorem \ref{thm: bsde existence} (iii) that
\[u(t,x)> \overline{u}(t,x), \quad \text{ for } (t,x) \in[0,T) \times (0,\infty)^d.\]

Before we state that both $u$ and $\overline{u}$ are viscosity solutions to \eqref{eq: pde}, we impose some additional assumptions. First, there exists a constant $\wt{K}$ such that
\begin{equation}\label{H: linear H}
 H(s, r) = \wt{K}(1+r) \quad \text{ for } (s,r) \in [0,T]\times (0,\infty). \tag{H6}
\end{equation}
This assumption implies that both $u(t,x)$ and $\overline{u}(t,x)$ are bounded from above by $C(1+\underline{x})$ on $[0,T] \times (0, \infty)^d$, where $C$ is a constant depending on $\mu, K, \wt{K}$, and $T$. Additionally,
\begin{equation}\label{H: nondeg}
 \sum_{i,j}^d (\sigma \sigma')_{ij}(x) v_i v_j >0, \quad \text{ for } x\in (0,\infty)^d \text{ and } v\in \Real^d\setminus \{0\}. \tag{H7}
\end{equation}
Denote
\[
\tau^x_n := \inf\{s\geq 0: X^x_s \notin \overline{\mathbb{B}_n}\} \wedge T.
\]
Since $X^x$ does not reach $\mathcal{O}$ in finite time, $X_{\tau_n^x} \in \mathbb{S}^+_n$ when $\tau^x_n <T$. Assumption \eqref{H: nondeg} implies that points on $\mathbb{S}^+_n$ are \emph{regular}, i.e., $\tau^x_n =0$ for any $x\in \mathbb{S}^+_n$ (see Theorem 2.3.3 in \cite{Pinsky}). This property will help us construct sequences of continuous functions which approximate $u$ and $\overline{u}$ from below.

Now we are ready to present the existence and uniqueness results for \eqref{eq: pde}.
\begin{thm}[Existence]\label{thm: pde existence}
 Suppose that \eqref{H: g linear growth} - \eqref{H: nondeg} hold. Then \eqref{eq: pde} admits two different viscosity solutions $u$ and $\overline{u}$. Both of them are nonnegative and bounded from above by $C(1+\underline{x})$, where $C$ depends on $\mu, K, \wt{K}$, and $T$. But $u(t,x) > \overline{u}(t,x)$ for $(t,x)\in [0,T)\times (0,\infty)^d$.
\end{thm}

\begin{rem}\label{rem: cont}
 Both $u$ and $\overline{u}$ are constructed via limits of increasing sequences of continuous functions. Therefore they are lower semi-continuous. When \eqref{eq: pde} is linear, the continuity of $u$ and $\overline{u}$ can be proved via the \emph{Schauder interior estimate} (see \cite{Ekstrom-Tysk}). When \eqref{eq: pde} is quasi-linear and the comparison result holds between viscosity super- and sub- solutions, $u=\overline{u}$ and they are continuous. A sufficient condition for the comparison result, hence the uniqueness result for \eqref{eq: pde}, is provided in Theorem \ref{thm: pde uniqueness} below.
\end{rem}

To obtain the comparison result for \eqref{eq: pde}, we need some additional assumptions:
for any $R>0$, there exists a function $m_R$ such that $\lim_{r\rightarrow 0} m_R(r) =0$ and
\begin{equation}\label{H: f cont x}
 |f(t,x,y,z) - f(t, x', y, z)| \leq m_R(|x-x'|(1+|z|)) \quad \text{ for } t\in [0,T], |x|, |x'|, |y| \leq R \text{ and } z\in \Real^d. \tag{H8}
\end{equation}
Additionally, we replace \eqref{H: f lip z} and \eqref{H: f mono y} with
\begin{align}
 & |f(t,x,y,z) - f(t,x,y,z')| \leq b(x)|z-z'|, \label{H: f lip z str}\tag{H3'.i}\\
 & |f(t,x,y,z) - f(t,x,y',z)| \leq \mu |y-y'|, \label{H: f lip y} \tag{H3'.ii}
\end{align}
for $t,x,y,y',z,z' \in [0,T]\times (0,\infty)^d \times \Real_+ \times \Real_+ \times \Real^d \times \Real^d$. Here $b(\cdot)$ is a bounded continuous function and $\mu$ is positive. We denote Assumptions \eqref{H: f lip z str}, \eqref{H: f lip y}, \eqref{H: f nonneg} and \eqref{H: f growth z} collectively as (H3').

As usual the uniqueness result follows from a comparison result. However, Theorem \ref{thm: pde existence} implies that the comparison result between viscosity super- and subsolutions fails when $X$ is a strict local martingale. To restore it, we assume the existence of a Lyapunov function $\Psi$, which ensures the martingale property of $X$.

\begin{thm}[Comparison]\label{thm: pde uniqueness}
 Suppose that (H1), (H2), (H3'), (H4) - (H8) hold. Moreover, there exists a strict positive function $\Psi: (0,\infty)^d \rightarrow (0,\infty)$ and a positive constant $\lambda$ such that
 \begin{align}
  &\cL \Psi(x) \leq \lambda (1+\Psi(x)), \quad \text{ for } x\in (0,\infty)^d, \label{H: Lpsi}\tag{H9.i}\\
  &\lim_{(0,\infty)^d \ni x\rightarrow \overline{x} \in \mathcal{O}} \Psi(x) = \infty, \label{H: psi close O}\tag{H9.ii}\\
  &\text{for any } M>0, \text{ there exists } R \text{ such that } \Psi(x) / \underline{x} \geq M \text{ for all } x \text{ with } \underline{x}\geq R, \label{H: Psi close inf}\tag{H9.iii}\\
  & c\Psi(x) \geq b(x) |\nabla \Psi(x) \sigma(x)|, \quad \text{ for some constant } c \text{ and all } x\in (0,\infty)^d. \label{H: Psi der bound}\tag{H9.iv}
 \end{align}
 Then for any nonnegative subsolution $u$ and supersolution $v$ which are both of at most linear growth in their spatial variables,
 \[
  u(t,x) \leq v(t,x), \quad \text{ for } (t,x)\in [0,T] \times (0,\infty)^d.
 \]
\end{thm}


\begin{rem}
 As we have seen in Remark \ref{rem: not reach O}, \eqref{H: Lpsi} and \eqref{H: psi close O} combined implies that $X$ never reaches $\mathcal{O}$ in finite time. On the other hand, \eqref{H: Lpsi} and \eqref{H: Psi close inf} ensure the martingale property of $X$. The reason is the following. \eqref{H: Lpsi} deduces that
 $
  \expec[\Psi(X^{t,x}_{s\wedge \tau_n})] \leq \Psi(x) + \lambda \int_0^s (1+\expec[\Psi(X^{t,x}_{u\wedge \tau_n})])\, du.
 $
 By Gronwall's inequality, the previous inequality yields $\expec[\Psi(X^{t,x}_{s\wedge \tau_n})] \leq (\Psi(x) + \lambda s) e^{\lambda s}=: M$ which is a constant independent of $n$. Now take any $\epsilon >0$, according to \eqref{H: Psi close inf}, there exists sufficiently large $R$ such that $\frac{\Psi(x)}{\underline{x}} \geq \frac{M}{\epsilon}$ for any $x$ such that $\underline{x}\geq R$. Then
 \[
  \expec\bra{\uX^{t,x}_{s\wedge \tau_n} \indic_{\{\uX^{t,x}_{s\wedge \tau_n} \geq R\}}} \leq \frac{\epsilon}{M} \expec\bra{\Psi(X^{t,x}_{s\wedge \tau_n}) \indic_{\{\uX^{t,x}_{s\wedge \tau_n} \geq R\}}} \leq M \frac{\epsilon}{M} =\epsilon, \quad \text{ for any } n.
 \]
 Hence $\{\uX^{t,x}_{s\wedge \tau_n}\}_{n\geq 0}$ is a uniformly integrable family. This implies that $\uX^{t,x}$, hence $X^{t,x}$, is a martingale.

 Assumption \eqref{H: Psi der bound} represents a balance between the growth restriction on $\sigma$ and the generator's dependence on $z$. Intuitively, the more restriction we put on the growth of $\sigma$, the wider class of generators Theorem \ref{thm: pde uniqueness} covers. Let us illustrate this point using the following examples.
\end{rem}

\begin{exa}[$\sigma$ has at most linear growth]\label{exa: at most linear growth}
 When $|\sigma(x)| \leq C(1+|x|)$ for some constant $C$ and all $x$, $b(\cdot)$ can be any bounded function, $\Psi$ can be chosen as $1+|x|^2$ (another function depending on the behavior of $\sigma$ near $\mathcal{O}$ needs to be added to $\Psi$ so that \eqref{H: psi close O} holds). One can check that \eqref{H: Lpsi}, \eqref{H: Psi close inf} and \eqref{H: Psi der bound} are satisfied. Therefore Theorem \ref{thm: pde uniqueness} holds for generators which are Lipschitz in $z$ and has bounded growth in $z$. In this case, the comparison result actually holds in the class of functions such that $\lim_{|x|\rightarrow \infty} |u(t,x)| e^{-A[\log |x|]^2} =0$ for some $A>0$ (see \cite{Barles-et-al}).
\end{exa}

\begin{exa}[No growth constraint on $\sigma$]\label{exa: mart only}
 If we know that $X$ is a martingale, but no other information on the growth of $\sigma$, Theorem \ref{thm: pde uniqueness} covers the case where the generator does not depend on $z$ (hence $b\equiv 0$). In fact, Assumption (H9) is sharp in 1 dimension: if $X$ is a 1-dimensional strict positive martingale, then there exists $\Psi$ which satisfies all conditions in (H9). Hence under other assumptions in Theorem \ref{thm: pde uniqueness}, the comparison holds among at most linear growth super- and subsolutions if and only if $X$ is a strict positive martingale.

 To construct $\Psi$, let us consider $\Psi_1(x) = 2\int_c^x dy \int_c^y \frac{dz}{\sigma^2(z)}$ for some $c>0$. It follows from Feller's test that $X$ does not reach $0$ in finite time if and only if $\lim_{x\downarrow 0} \Psi_1(x) =\infty$. On the other hand, $X$ is a martingale if and only if $\int_c^\infty \frac{x}{\sigma^2(x)} dx =\infty$ (see \cite{Delbaen-Shirakawa}). Then consider $\Psi_2(x) = x+ \int_c^x dy \int_c^y \frac{z}{\sigma^2(z)} dz$. We set $\Psi = \Psi_1 + \Psi_2$. \eqref{H: psi close O} and \eqref{H: Psi der bound} clearly hold; \eqref{H: Psi close inf} follows from the fact that $\lim_{x\rightarrow \infty} \Psi_2(x) / x = \lim_{x\rightarrow \infty} \Psi_2'(x) =\infty$; \eqref{H: Lpsi} follows from $\cL \Psi_1 =1$ and $\cL \Psi_2 \leq \Psi_2/2$.
\end{exa}

\begin{exa}[$\sigma$ has superlinear growth but $X$ is still a martingale]\label{exa: super linear growth}
 Consider the 1-dimensional SDE $dX_t = \sigma(X_t) dB_t$ where $\sigma(x) = \left\{\begin{array}{ll}x & \text{ if } x\leq\textit{e} \\ x\sqrt{\log x} & \text{ if } x> \textit{e} \end{array}\right.$. One can check that $\sigma$ is locally Lipschitz in $(0,\infty)$ and the solution $X$ does not reach $0$ in finite time, because $X$ is a Geometric Brownian motion when $X\leq \textit{e}$. On the other hand, since $\int_{\textit{e}}^\infty \frac{x}{x^2 \log x} dx =\infty$, $X$ is a martingale (see \cite{Delbaen-Shirakawa}). Consider $b(x) = \left\{\begin{array}{ll} 1 & \text{ if } x\leq \textit{e} \\ \frac{\textit{e}}{x \sqrt{\log{x}}} & \text{ if } x>\textit{e} \end{array}\right.$. We will show in the next paragraph that $\Psi$ exists and (H9) is satisfied. Then Theorem \ref{thm: pde uniqueness} holds in this case, where the generator may depend on $z$ nontrivially.

 Let us set
 \[
  \Psi(x) = \frac{1}{x} + x+ \int_{\textit{e}}^x dy \int_{\textit{e}}^y \frac{z}{\sigma^2(z)} dz.
 \]
 Clearly \eqref{H: psi close O} holds and so does \eqref{H: Psi close inf}, which follows from the same argument as in the last example. Now we are going to verify \eqref{H: Lpsi} and \eqref{H: Psi der bound}. First,
 \[
  \frac12 \sigma^2(x) \Psi^{''}(x) = \left\{\begin{array}{ll} \frac1x + \frac12 x & \text{ if } x\leq \textit{e}\\ \frac{\log{x}}{x} + \frac12 x & \text{ if } x> \textit{e}\end{array}\right. \leq \Psi(x).
 \]
 Then \eqref{H: Lpsi} holds. Second,
 \[
  b(x)|\Psi'(x)|\sigma(x) \leq \left\{\begin{array}{ll}\frac1x + 2x - x\log{x}& \text{ if } x\leq \textit{e} \\ e\pare{\frac{1}{x^2} + 1+ \log\log x} & \text{ if } x> \textit{e}\end{array}\right. \leq C \Psi(x),
 \]
 where the second inequality holds for sufficiently large $C$ because $\lim_{x\downarrow 0} x\log{x} =0$ and $\lim_{x\rightarrow \infty} \frac{\int_{\textit{e}}^x \log\log y dy}{\log\log x} = \infty$ from l'Hopital rule. Hence \eqref{H: Psi der bound} is also verified.
\end{exa}

\section{Construction of multiple solutions to \eqref{eq: bsde}}\label{sec: bsde solution}
Let us first discuss the construction of $(Y, Z)$ and $(\oY, \oZ)$ intuitively. Recall
$\tau_n = \inf\{s\geq 0 \such X_s \notin \overline{\mathbb{B}_n}\}\wedge T$ for $n>0$. The supermartingale property of $X$ implies that $\{\tau_n = T\}$ increases to $\Omega$ as $n\rightarrow \infty$. Moreover, the stopped processes $X_{\cdot \wedge \tau_n}$ and $\uX_{\cdot \wedge \tau_n}$ are martingales. Given a sequence of random variables $\xi_n \in \F_{\tau_n}$, we consider the following sequence of BSDEs:
\begin{equation}\label{eq: bsde Yn}
 Y^n_t = \xi_n + \int_t^T \indic_{\{s\leq \tau_n\}} f(s, X_s, Y^n_s, Z^n_s) \,ds - \int_t^T Z_s^n\, dB_s, \quad \text{ for each } n\geq 0.
\end{equation}
To approximate \eqref{eq: bsde}, we choose two different sequences of terminal conditions for the previous BSDE:
\begin{equation}\label{eq: def xi oxi}
 \xi_n := g(X_{\tau_n}) \quad \text{ and } \quad \overline{\xi}_n := g_n(X_{\tau_n}),
\end{equation}
where $g_n(x):= g(x) h_n(x)$ and $h_n(\cdot)$ is a continuous function such that $0\leq h_n \leq 1$ and $h_n(x) = \left\{\begin{array}{ll}1, & x\in \mathbb{B}_{n-1} \\ 0, & x\notin \mathbb{B}_n \end{array}\right.$. Since $g$ is bounded on $\overline{\mathbb{B}_n}$, both $\xi_n$ and $\overline{\xi}_n$ are bounded. Then under Assumptions \eqref{H: f lip z} and \eqref{H: f mono y}, \eqref{eq: bsde Yn} admits a solution: $(Y^n, Z^n)$ when the terminal condition is $\xi_n$; $(\oY^n, \oZ^n)$ when the terminal condition is $\overline{\xi}_n$. Both these solutions are also unique inside the class $(\cS^\infty, \cM^2)$. See e.g. Theorem 2.2 and Proposition 2.2 in \cite{Pardoux}.

Notice that $\overline{\xi}_n = g_n(X_T) \indic_{\{\tau_n =T\}}$. Both $\{\xi_n\}_{n\geq 0}$ and $\{\overline{\xi}_n\}_{n\geq 0}$ converge to $g(X_T)$ in probability as $n\rightarrow \infty$. This convergence motivates us to construct $Y$ and $\oY$ via limits of $\{Y^n\}_{n\geq 0}$ and $\{\oY^n\}_{n\geq 0}$, respectively. It is important to note that the convergence of $\{\xi_n\}_{n\geq 0}$ and $\{\overline{\xi}_n\}_{n\geq 0}$ is in probability, not necessarily in expectation. This allows that $\{Y^n\}_{n\geq 0}$ and $\{\oY^n\}_{n\geq 0}$ eventually converge to different solutions. To make this idea rigorous, we will employ a localization argument in \cite{Briand-Hu} and then apply the monotone stability result for solutions of BSDE in \cite{Kobylanski}. Before carrying out these steps, let us prepare the following two lemmas.

\begin{lem}\label{lem: increasing}
Let (H1), (H2), (H3.i) - (H3.iii) hold. Then
\[
 \oY^{n+1}_t \geq \oY^n_t, \quad \text{ for } t\in [0,T].
\]
If (H4.i) also holds,
\[
 Y^{n+1}_t \geq Y^n_t,  \quad \text{ for } t\in [0, \tau_n].
\]
\end{lem}

\begin{proof}
Recall $\overline{\xi}_n = g_n(X_T) \indic_{\{\tau_n =T\}}$. Since both $\{g_n\}_{n\geq 0}$ and $\{\tau_n\}_{n\geq 0}$ are nondecreasing, then $\{\overline{\xi}_n\}_{n\geq 0}$ is also nondecreasing. On the other hand, $\indic_{\{s\leq \tau_n\}} f \leq \indic_{\{s\leq \tau_{n+1}\}} f$ since $f$ is nonnegative. Therefore the first statement follows from the comparison theorem (see e.g. Theorem 2.4 in \cite{Pardoux}) directly.

 To prove the second statement, we first show
 \begin{equation}\label{eq: xi-inc}
  \expec\bra{\xi_{n+1} \,|\, \F_{\tau_n}} \geq \xi_n.
 \end{equation}
 Indeed, this follows from
 \[
  \expec\bra{\xi_{n+1} \,|\, \F_{\tau_n}} = K \pare{1+ \expec[\uX_{\tau_{n+1}} \,|\, \F_{\tau_n}]} - \expec\bra{\og(X_{\tau_{n+1}}) \,|\, \F_{\tau_n}} \geq K(1+\uX_{\tau_n}) -\og(X_{\tau_n}) = \xi_n,
 \]
 where the inequality uses \eqref{H: og supermart} and the martingale property of $\underline{X}_{\tau_{n+1} \wedge \cdot}$. Now consider the following BSDE:
 \begin{equation}\label{eq: tY}
  \wt{Y}^{n+1}_t = \xi_{n+1} + \int_t^T \indic_{\set{s\leq \tau_n}} f(s, X_s, \wt{Y}^{n+1}_s, \wt{Z}^{n+1}_s) \,ds - \int_t^T \wt{Z}^{n+1}_s \, dB_s.
 \end{equation}
 It admits a unique solution $(\wt{Y}^{n+1}, \wt{Z}^{n+1})\in (\cS^\infty, M^2)$. Since $\tau_n \leq \tau_{n+1}$ and $f\geq 0$, the comparison theorem implies that
 \begin{equation}\label{eq: Y>tY}
  Y^{n+1}_t \geq \wt{Y}^{n+1}_t, \quad \text{ for } t\in [0,T].
 \end{equation}
 Taking conditional expectation with respect to $\F_{\tau_n}$ on both sides of \eqref{eq: tY} gives
 \[
  \wt{Y}^{n+1}_t = \expec[\xi_{n+1} \,|\, \F_{\tau_n}] + \int_t^T \indic_{\set{s\leq \tau_n}} f(s, X_s, \wt{Y}^{n+1}_s, \wt{Z}^{n+1}_s)\, ds - \int_t^T \wt{Z}^{n+1}_s \indic_{\set{s\leq \tau_n}} dB_s, \quad t\in [0, \tau_n].
 \]
 Compare the previous BSDE with the one satisfied by $Y^n$. The comparison theorem and \eqref{eq: xi-inc} combined gives
 \[
  \wt{Y}^{n+1}_t \geq Y^n_t, \quad \text{ for } t\in [0,\tau_n].
 \]
 Then the second statement follows after combining the previous inequality with \eqref{eq: Y>tY}.
\end{proof}

The following lemma gives a upper bound for $Y^n$ and $\oY^n$.

\begin{lem}\label{lem: Y bound}
 Let \eqref{H: f mono y} and \eqref{H: f growth z} hold. For any $n\geq 0$,
 \[
  Y^n_t\leq C \pare{\expec\bra{\left. \xi_n + \int_t^T H(s, \uX_s) \,ds \right|\F_t}}\leq C \pare{\expec[\xi_n \,|\, \F_t] + \int_t^T H(s, \uX_t) ds}, \quad t\in [0,T],
 \]
 where $C = e^{(\mu \vee 0) T}$. The same statement holds for $(\oY^n, \overline{\xi}_n)$ as well.
\end{lem}

\begin{proof}
 We only prove the statement for $(Y^n, \xi_n)$, the same argument applies to the statement for $(\oY^n, \overline{\xi}_n)$ as well.
 Consider the following ODE:
 \[
  \varphi^n_t = \xi_n + \int_t^T \indic_{\set{s\leq \tau_n}} (H(s, \uX_s) + \mu \varphi^n_s) \, ds
 \]
 and define $\Phi^n_t := \expec[\varphi^n_t \,|\, \F_t]$. The solution to the previous ODE is
 \[
  \varphi^n_t = \xi_n \quad \text{for } t\geq \tau_n \quad \text{ and } \quad \varphi^n_t = e^{\mu(\tau_n -t)}\xi_n + \int_t^{\tau_n} e^{\mu(s-t)} H(s, \uX_s) ds, \quad \text{for } t<\tau_n.
 \]
 It then follows $0\leq \varphi^n_t \leq C(\xi_n + \int_t^T H(s, \uX_s)) ds$, which yields
 \[
  0\leq \Phi^n_t \leq C\pare{\expec\bra{\left.\xi_n + \int_t^T H(s, \uX_s) \, ds \right| \F_t}}.
 \]
 Since $r \mapsto H(\cdot, r)$ is concave and nondecreasing,
 \[
  \expec\bra{\left.\int_t^T H(s, \uX_s) \,ds\right|\F_t} = \int_t^T \expec\bra{\left.H(s, \uX_s) \right| \F_t} ds \leq \int_t^T H\pare{s, \expec[\uX_s \,|\, \F_t]} \, ds \leq \int_t^T H(s, \uX_t) \,ds,
 \]
 where the second inequality follows from the supermartingale property of $\uX$.
 Therefore the last two estimates combined gives
 \[
  0\leq \Phi^n_t \leq C\pare{\expec[\xi_n \,|\, \F_t] + \int_t^T H(s,\uX_t) \,ds}.
 \]
 Now the statement follows if we can show
 \[
  \Phi^n_t \geq Y^n_t, \quad t\in [0,T].
 \]
 To this end, note that \eqref{H: f mono y} and \eqref{H: f growth z} imply $f(t,x,y,z) \leq H(t, \underline{x}) + \mu y$ for any $(t,x,y,z)\in [0,T]\times \Real_+^d \times \Real_+ \times \Real^d$. Then the previous claim follows from the same comparison argument in Lemma 1 of \cite{Briand-Hu}.
\end{proof}

Now we are ready to prove the first main result.
\begin{proof}[Proof of Theorem \ref{thm: bsde existence}]
 The proof is split into several steps.

 \underline{\textit{Step 1: Construction of solutions.}}
 We will only present the construction of $(Y, Z)$ from the limit of $\{(Y^n, Z^n)\}_{n\geq 0}$. The solution $(\oY, \oZ)$ can be similarly constructed via the limit of $\{(\oY^n, \oZ^n)_{n\geq 0}\}$. Fix $k\in \Natural$. We stop every $(Y^n, Z^n)$ at $\tau_k$ by defining
 \[
  Y^{n,k}_t := Y^n_{t\wedge \tau_k} \quad \text{and} \quad Z^{n,k}_t := Z^n_t \indic_{\{t\leq \tau_k\}}.
 \]
 These stopped processes satisfy the following BSDE:
 \[
  Y^{n,k}_t = Y^n_{\tau_k} + \int_t^T \indic_{\{s\leq \tau_k\}} f(s, X_s, Y^{n,k}_s, Z^{n,k}_s) \,ds -\int_t^T Z^{n,k}_s \, dB_s.
 \]
 It follows from Lemma \ref{lem: increasing} that $\{Y^{n,k}\}_{n\geq k}$ is a nondecreasing sequence. Moreover $\{Y^{n,k}\}_{n\geq k}$ is bounded uniformly in $n$. Indeed, Lemma \ref{lem: Y bound} and \eqref{H: g linear growth} implies that
 \[
 \begin{split}
  0\leq Y^{n,k}_t &= Y^n_{t\wedge \tau_k} \leq C\pare{K(1+ \uX_{t\wedge \tau_k}) + \int_{t\wedge \tau_k}^T H(s, \uX_{t\wedge \tau_k}) \,ds} \leq M_k, \quad t\in[0,T].
 \end{split}
 \]
 Here $M_k$, depending on the maximum of $H$ on $[0,T] \times [\min_{x\in \overline{\mathbb{B}_k}} \underline{x}, \max_{x\in \overline{\mathbb{B}_k}} \underline{x}]$, is a constant independent of $n$.

 Since $\{Y^{n,k}\}_{n\geq k}$ is monotone and uniformly bounded, it follows from Proposition 2.4 in \cite{Kobylanski} that $\{Y^{n,k}\}_{n\geq k}$ converges uniformly on $[0,T]$ to a continuous process $\kY_\cdot:= \lim_{n\rightarrow \infty} Y^{n,k}_\cdot$ and $\{Z^{n,k}\}_{n\geq k}$ converges to some $\kZ$ in $\cM^2$, such that $(\kY, \kZ)\in (\cS^\infty, \cM^2)$ is a solution to the following BSDE:
 \begin{equation}\label{eq: bsde-kY}
  \kY_t = \eta_k +\int_t^T \indic_{\set{s\leq \tau_k}} f(s, X_s, \kY_s, \kZ_s)\, ds -\int_t^T \kZ_s\, dB_s,
 \end{equation}
 where $\eta_k =\lim_{n\rightarrow \infty} Y^n_{\tau_k}$. Note $\eta_k=g(X_T)$ when $\tau_k=T$. We will use this observation later.

 Now coming back to the definition of $Y^{n,k}$ and $\kY$, we have
 \[
  {}^{k+1}\kern-0.15em Y_{t\wedge \tau_k} =\lim_{n\rightarrow \infty} Y^{n, k+1}_{t\wedge \tau_k} = \lim_{n\rightarrow \infty} Y^n_{t\wedge \tau_k \wedge \tau_{k+1}} = \lim_{n\rightarrow \infty} Y^n_{t\wedge \tau_k} =\lim_{n\rightarrow \infty} Y^{n,k}_t = \kY_t.
 \]
 On the other hand, it follows from $\lim_{n\rightarrow \infty} \expec[\int_0^T |{}^{k}\kern-0.15em Z_s - Z^{n,k}_s|^2 \,ds] =0$ that $\lim_{n\rightarrow \infty} \expec[\int_0^{\tau_k} |{}^{k}\kern-0.15em Z_s - Z^{n,k}_s|^2 \,ds]=0$. Similarly, $\lim_{n\rightarrow \infty} \expec[\int_0^{\tau_k} |{}^{k+1}\kern-0.15em Z_s - Z^{n,k+1}_s|^2 \,ds]=0$. Noticing that $Z^{n,k}_s \indic_{\{s\leq \tau_k\}} = Z^n_s \indic_{\{s\leq \tau_k\}} = Z^{n, k+1}_s \indic_{\{s\leq \tau_k\}}$, we obtain $\expec[\int_0^{\tau_k} |{}^{k+1}\kern-0.15em Z_s - {}^{k}\kern-0.15em Z_s|^2 \,ds] =0$.
 Therefore we can define $Y$ and $Z$ via
 \[
  Y_{t\wedge \tau_k} := \kY_t \quad \text{ and } \quad Z_t := \kZ_t \quad \text{if } t\in [0,\tau_k].
 \]
 When $\tau_k=T$, since $\kY$ is continuous on $[0,T]$, so is $Y$. Moreover $\lim_{t\rightarrow T} Y_t = \lim_{t\rightarrow T} \kY_t = \eta_k = g(X_T)$ on $\set{\tau_k=T}$. By sending $k$ to infinity and recalling that $\bigcup_{k\in \Natural} \{\tau_k = T\} = \Omega$, we deduce that $Y$ is almost surely continuous and  $\lim_{t\rightarrow T} Y_t = g(X_T)$. On the other hand, from the definition of $Z$,
 \[
 \begin{split}
  \prob\pare{\int_0^T |Z_s|^2 \,ds =\infty} &= \prob\pare{\int_0^T |Z_s|^2 \,ds =\infty, \tau_k=T} + \prob\pare{\int_0^T |Z_s|^2 ds =\infty, \tau_k<T} \\
  &\leq \prob\pare{\int_0^{\tau_k} |\kZ_s|^2 \,ds = \infty} + \prob(\tau_k<T).
 \end{split}
 \]
 The right-hand-side of the previous inequality converges to zero as $k\rightarrow \infty$. Therefore $\int_0^T Z_s^2 \,ds <\infty$. Now following from \eqref{eq: bsde-kY}, $(Y,Z)$ satisfies
 \[
  Y_{t\wedge \tau_k} = Y_{\tau_k} +\int_{t\wedge \tau_k}^{\tau_k} f(s, X_s, Y_s, Z_s) \,ds -\int_{t\wedge \tau_k}^{\tau_k} Z_s dB_s.
 \]
 Sending $k$ to infinity, we conclude that $(Y,Z)$ is a solution to \eqref{eq: bsde}.

 \vs
 \underline{\textit{Step 2: Uniform integrability.}} From Lemma \ref{lem: Y bound},
 \begin{equation}\label{eq: ub-oY}
 \begin{split}
  \oY_t &= \lim_{n\rightarrow \infty} \oY^n_t \leq C\pare{\lim_{n\rightarrow \infty} \expec[g_n(X_T)\indic_{\{\tau_n =T\}} \,|\, \F_t] + \expec\bra{\left.\int_t^T H(s, \uX_s) \,ds\right| \F_t}} \\
  &= C\pare{\expec\bra{\left.g(X_T) + \int_t^T H(s, \uX_s) \,ds\right| \F_t}}, \quad \text{ on } \{t\leq \tau_k\},
 \end{split}
 \end{equation}
 where the second equality follows from the dominated convergence theorem. Send $k$ to infinity, \eqref{eq: ub-oY} holds for $t\in [0,T]$. Therefore, $\oY$ is of class \emph{D} because it is nonnegative and bounded from above by a uniformly integrable martingale. On the other hand, combined with \eqref{H: g linear growth}, \eqref{eq: ub-oY} also implies $\oY\in \C$.

 Now let us switch our attention to $Y$. First, since $Z^n \in \cM^2$, $\int_0^\cdot Z_s^n \,dB_s$ is a martingale. Then $f\geq 0$ implies that $Y^n_t =\expec\bra{\xi_n + \int_t^T \indic_{\{s\leq \tau_n\}} f(s, X_s, Y^n_s, Z^n_s) ds \,|\, \F_t} \geq \expec[\xi_n\,|\,\F_t]$. The construction of $Y$ then yields
 \[
  Y_t = \lim_{n\rightarrow \infty} Y_t^n \geq \lim_{n\rightarrow \infty} \expec[g(X_{\tau_n}) \,|\, \F_t].
 \]
 In order to derive $\lim_{n\rightarrow \infty}\expec\bra{g(X_{\tau_n}) \,|\, \F_t}$, recall $\overline{g}(x) = K(1+\underline{x})- g(x)$ and $\overline{g}(x) \leq \overline{G}(\underline{x})$ from \eqref{H: og strict sublinear}. Since $\overline{G}$ is nondecreasing and $\lim_{r\rightarrow \infty} \overline{G}(r)/r =0$, there exists a function $\psi$ such that $\psi(\overline{G}(r))\leq r$ for $r\geq 0$ and $\lim_{y\uparrow \infty} \psi(y)/y =\infty$. Indeed, set $\psi(y) = \inf\{r\geq 0: \overline{G}(r) \geq y\}$. $\psi$ is nondecreasing, $\lim_{y\uparrow \infty} \psi(y) = \infty$, and $\psi(\overline{G}(r))\leq r$. On the other hand, since $\overline{G}(\psi(y)) \geq y$, it follows $0= \lim_{y\uparrow \infty} \frac{\overline{G}(\psi(y))}{\psi(y)} \geq \limsup_{y\uparrow \infty} \frac{y}{\psi(y)}$. Therefore $\lim_{y\uparrow \infty} \psi(y) / y = \infty$. These properties on $\psi$ imply that $\expec[\psi(\overline{G}(\uX_{\tau_n}))] \leq \expec[\uX_{\tau_n}] \leq \underline{x}$ for any $n$. From de la Vall\'{e}e Poussin criteria (see Lemma 3 in \cite{Shiryaev} pp. 190), the previous inequalities imply that $\{\overline{G}(\uX_{\tau_n})\}_{n\geq 0}$, hence $\{\overline{g}(X_{\tau_n})\}_{n\geq 0}$, is uniformly integrable. As a result, $\lim_{n\rightarrow \infty} \expec[\overline{g}(X_{\tau_n}) \,|\, \F_t] = \expec[\overline{g}(X_T) \,|\, \F_t]$. Go back to the limit of $\expec[g(X_{\tau_n})\,|\, \F_t]$,
 \begin{equation}\label{eq: limit xi}
 \begin{split}
  \lim_{n\rightarrow \infty} \expec\bra{g(X_{\tau_n})\,|\, \F_t} &= K\pare{1+ \lim_{n\rightarrow \infty} \expec\bra{\uX_{\tau_n} \,|\, \F_t}} - \lim_{n\rightarrow \infty} \expec\bra{\overline{g}(X_{\tau_n}) \,|\, \F_t} \\
  &= K\pare{1+ \lim_{n\rightarrow \infty} \uX_{t\wedge \tau_n}} - \expec\bra{\overline{g}(X_T) \,|\, \F_t} \\
  & = K(1+\uX_t) - \expec\bra{\overline{g}(X_T) \,|\, \F_t} \\
  &= K\pare{\uX_t - \expec[\uX_T \,|\, \F_t]} + \expec\bra{g(X_T) \,|\, \F_t}.
 \end{split}
 \end{equation}
 Now if $Y$ was of class \emph{D}, $\uX$ would also be, since $\expec\bra{g(X_T) - K \uX_T \,|\, \F_\cdot}$ is already of class \emph{D}. However this contradicts with the strict local martingale property of $X$.

 Applying the similar estimate to the upper bound of $Y$, we obtain
 \begin{equation}\label{eq: ub-Y}
 \begin{split}
  Y_t &\leq C\pare{K\pare{\uX_t - \expec[\uX_T \,|\, \F_t]} + \expec\bra{\left.g(X_T) + \int_t^T H(s, \uX_s) \,ds\right| \F_t}}\\
  & \leq C\pare{K(1+ \uX_t) + \expec\bra{\left.\int_t^T H(s, \uX_s)\,ds\right| \F_t}},
 \end{split}
 \end{equation}
 where the second inequality holds since $g(x)-K \underline{x} \leq K$. Therefore $Y\in \C$ follows from the previous inequality.

 \underline{\textit{Step 3: $(\oY, \oZ)$ is the minimal solution.}} Since $\oZ^n \in \cM^2$, it follows from the definition of $\oY^n$ that $\oY^n_{\tau_n} =\expec\bra{g_n(X_T) \indic_{\{\tau_n =T\}} \,|\, \F_{\tau_n}}$. On the other hand, for any solution $(\wt{Y}, \wt{Z})$ to \eqref{eq: bsde} such that $\wt{Y}\in \C$,
 \[
 \wt{Y}_{t\wedge \zeta_n} = \wt{Y}_{u \wedge \zeta_n} + \int_t^u \indic_{\{s\leq \zeta_n\}} f(s, X_s, \wt{Y}_s, \wt{Z}_s) \,ds - \int_t^u \wt{Z}_s \indic_{\{s\leq \zeta_n\}}\,dB_s, \quad 0\leq t\leq u \leq T,
 \]
 where $\zeta_n$ is chosen as $\inf\{t\geq 0\such |\int_0^t |\wt{Z}_u|^2 du \geq n\}$. Since $f$ is nonnegative, the previous BSDE gives $\wt{Y}_t \indic_{\{t\leq \zeta_n\}} \geq \expec[\wt{Y}_{u\wedge \zeta_n} \,|\, \F_t] \indic_{\{t\leq \zeta_n\}}$. Sending $n\rightarrow \infty$ and utilizing Fatou's lemma, we obtain $\wt{Y}_t \geq \expec[\wt{Y}_u \,|\, \F_t]$, hence $\wt{Y}$ is a supermartingale. As a result,
 \begin{equation}\label{eq: tY>oY}
  \wt{Y}_{\tau_n} \geq \expec[\wt{Y}_T \,|\, \F_{\tau_n}] = \expec[g(X_T) \,|\, \F_{\tau_n}] \geq \expec[g_n(X_T) \indic_{\{\tau_n =T\}} \,|\, \F_{\tau_n}] = \oY^n_{\tau_n}.
 \end{equation}
 On the other hand, since $\wt{Y}\in \C$ and $r\mapsto H(\cdot, r)$ is concave and nondecreasing,
 \[
  \wt{Y}_t \leq C\pare{K(1+ \uX_t) + \int_t^T H(s, \uX_t) \,ds}.
 \]
 Therefore $\wt{Y}_{\cdot \wedge \tau_n} \in \cS^\infty$, which implies $\wt{Z}_\cdot \indic_{\{\cdot \leq \tau_n\}} \in \cM^2$ (see Proposition 2.2 in \cite{Pardoux}).

 Now compare the following two BSDEs:
 \begin{align*}
  \oY^n_t & = \oY^n_{\tau_n} +\int_t^T \indic_{\{s\leq \tau_n\}} f(s, X_s, \oY^n_s, \oZ^n_s) \, ds - \int_t^T \oZ^n_s \indic_{\{s\leq \tau_n\}} dB_s,\\
  \wt{Y}_t & = \wt{Y}_{\tau_n} + \int_t^T \indic_{\{s\leq \tau_n\}} f(s, X_s, \wt{Y}_s, \wt{Z}_s)\, ds - \int_t^T \wt{Z}_s \indic_{\{s\leq \tau_n\}} dB_s.
 \end{align*}
 Thanks to \eqref{eq: tY>oY}, the comparison theorem in $(\cS^\infty, \cM^2)$ (see e.g. Theorem 2.4 in \cite{Pardoux}) implies that
 \[
  \wt{Y}_t \geq \oY^n_t, \quad t\in [0, \tau_n].
 \]
 Since the choice of $n$ is arbitrary, $\wt{Y}_t \geq \oY_t$ for $t\in [0,T]$ is then clear.

 \underline{\textit{Step 4: $Y_0 > \oY_0$.}} Let us define
 \[
  \alpha_t = \left\{\begin{array}{ll}(Y_t - \oY_t)^{-1} \pare{f(t, X_t, Y_t, \oZ_t) - f(t, X_t, \oY_t, \oZ_t)} & \text{ if } Y_t \neq \oY_t\\ 0 & \text{ if } Y_t = \oY_t\end{array}\right.,
 \]
 and the $\Real^d$- valued process $\{\beta_{t}; t\in [0,T]\}$ as follows. For $1\leq i\leq d$, let $Z_t^{(i)}$ be the $d$-dimensional vector whose first $i$ components are equal to those of $\oZ_t$ and whose last $d-i$ components are equal to those of $Z_t$. Then we define for $1\leq i\leq d$,
 \[
  \beta^i_t = \left\{\begin{array}{ll}(\oZ^i_t-Z_t^i) \pare{f(t, X_t, Y_t, \oZ^{(i)}_t) - f(t, X_t, Y_t, Z^{(i-1)}_t)} & \text{ if } Z_t^i \neq \oZ^i_t\\ 0 & \text{ if } Z^i_t = \oZ^i_t\end{array}\right..
 \]
 Note that $\{\alpha_t; t\in [0,T]\}$ and $\{\beta_t; t\in [0,T]\}$ are both progressively measurable, $\alpha_t \leq \mu$, and $|\beta_t| \leq \nu$ from \eqref{H: f lip z} and \eqref{H: f mono y}.

 For $0\leq t\leq u\leq T$, define $\Gamma_{t, u} := \exp\pare{\int_t^u \pare{\alpha_s - \frac12 |\beta_s|^2} ds + \int_t^u \beta'_s dB_s}$. Then $(\cY_t, \cZ_t) := (Y_t - \oY_t, Z_t - \oZ_t)$ satisfies
 \begin{equation}\label{eq: bsde Y-oY}
  \cY_t = \Gamma_{t, u}\cY_u -\int_t^u \Gamma_{t,s}(\cZ_s + \cY_s \beta'_s) dB_s.
 \end{equation}
  Set $\zeta_n = \inf\{t\geq 0 \such \int_0^t\Gamma_{0,u}^2|\cZ_u+ \cY_u \beta'_u|^2 du \geq n\}$. We have from \eqref{eq: bsde Y-oY} that
 \[
  \cY_t \indic_{\{t\leq \zeta_n\}} = \expec[\Gamma_{t, u\wedge \zeta_n} \cY_{u\wedge \zeta_n} \,|\, \F_t] \indic_{\{t\leq \zeta_n\}}.
 \]
 Since both $Y$ and $\oY$ are continuous processes, moreovoer $Y_t \geq \oY_t \geq 0$ for $t\in [0,T]$ from Step 3, we have $\prob(\cY_t \geq 0 \text{ for all } t\in [0,T])=1$. Therefore $\cY_{\cdot \wedge \zeta_n} \geq 0$ for any $n\geq 0$. It then follows from Fatou's lemma that
 \[
  \cY_t = \lim_{n\rightarrow \infty}\expec[\Gamma_{t, u\wedge \zeta_n} \cY_{u\wedge \zeta_n} \,|\, \F_t] \indic_{\{t\leq \zeta_n\}} \geq \expec[\Gamma_{t, u} \cY_u \,|\, \F_t].
 \]
 As a result, $\{\Gamma_{0,t}\cY_t; t\in [0,T]\}$ is a nonnegative super martingale.

 Now if $\cY_0 = Y_0-\oY_0 =0$, then $\cY_t =0$ for any $t\in [0,T]$, which implies that $\prob(\cY_t = 0 \text{ for all } t\in [0,T])=1$. However, this contradicts with the fact that $\oY$ is of class $D$ but $Y$ is not.
\end{proof}

\begin{proof}[Proof of Corollary~\ref{cor: inf soln}]
 Consider $\{\xi^\alpha_n\}_{\alpha \in [0,1], n\geq 1}$ where $\xi^\alpha_n = (1-\alpha) g_n(X_{\tau_n}) + \alpha g(X_{\tau_n})$. Since $g_n(X_{\tau_n}) = g_n(X_T) \indic_{\{\tau_n = T\}} \leq g(X_{\tau_n})$, then $\{\xi_n^\alpha\}_{\alpha \in [0,1]}$ is nondecreasing. It is also clear that $\xi_n^0 = \overline{\xi}_n$ and $\xi^1_n = \xi_n$. Consider \eqref{eq: bsde Yn} whose terminal condition is replaced by $\xi_n^\alpha$. We denote its solution by $(Y^{n,\alpha}, Z^{n,\alpha})$. Walking through Lemma \ref{lem: increasing} and Step 1 in Theorem \ref{thm: bsde existence}, we obtain a sequence of \eqref{eq: bsde} solutions $(Y^\alpha, Z^\alpha)_{\alpha \in [0,1]}$ such that $\{Y^\alpha\}_{\alpha\in[0,1]}$ is nondecreasing. it is also clear that $(Y^0, Z^0) = (\overline{Y}, \overline{Z})$ and $(Y^1, Z^1) = (Y, Z)$.

 In this paragraph, we will show $\{Y^\alpha_0\}_{\alpha\in [0,1]}$ is strictly increasing. For any $0\leq \alpha < \alpha' \leq 1$, applying the argument in Step 4 of Theorem \ref{thm: bsde existence} to $Y^{n, \alpha}$ and $Y^{n, \alpha'}$, we obtain
 \begin{equation}\label{eq: Y^n est}
  \mathcal{Y}^n_0 \geq \expec\bra{\Gamma_{0, \tau_n} \mathcal{Y}^n_{\tau_n}} \geq C_\mu \expec[\mathcal{Y}^n_{\tau_n}] = C_\mu (\alpha'-\alpha)\expec[g(X_{\tau_n}) - g_n(X_{\tau_n})],
 \end{equation}
 where $\mathcal{Y}^n = Y^{n, \alpha'} - Y^{n, \alpha}$ and $\Gamma_{0, \tau_n} \geq C_\mu =: \exp(-\mu T)$ because $f(t,x,y,z)$ does not depend on $z$ and is Lipschitz in $y$ with some Lipschitz constant $\mu$. Sending $n\rightarrow \infty$ in \eqref{eq: Y^n est}, and utilizing arguments in Step 2 of Theorem \ref{thm: bsde existence}, we obtain
 \[
  Y^{\alpha'}_0 - Y^\alpha_0 \geq C_\mu (\alpha'-\alpha) K(\underline{X}_0 - \expec[\underline{X}_T])>0,
 \]
 since $X$ is a nonnegative strict local martingale. This confirms the claim.
\end{proof}

In the rest of this section, we will prove that any solution, whose first component is in $\mathcal{C}$, is inside the class $(\cS^p, \cM^p)$ for any $p\in (0,1)$. Let us first recall the following version of Doob's inequality.

\begin{lem}\label{lem: local mart Lp}
 If $L$ is a 1-dimensional nonnegative local martingale, then
 \[
  \expec[\sup_{0\leq t\leq T} L_t^p] \leq \frac{1}{1-p} L_0^p, \quad \text{for any } p\in (0,1).
 \]
\end{lem}
\begin{proof}
 Being a nonnegative local martingale implies that $L$ is a supermartingale. It then follows from Doob's second submartingale inequality (see Theorem 1.3.8 in \cite{Karatzas-Shreve-BM}) that
 \[
  y \prob(\inf_{0\leq t\leq T} (-L_t) \leq -y) \leq \expec[(-L_T)^+] - \expec[-L_0] = L_0.
 \]
 Denote $L_* = \sup_{0\leq t\leq T} L_t$, it then follows
 \[
  \expec[L_*^p] = \expec\bra{\int_0^\infty \indic_{\set{L_* \geq y}} p y^{p-1} \, dy} \leq \int_0^\infty \min\left\{1, L_0/y\right\} p y^{p-1} \, dy = \frac{L_0^p}{1-p}.
 \]
\end{proof}

\begin{prop}\label{prop: integrability}
 For any solution $(Y, Z)$ to \eqref{eq: bsde} such that $Y\in \C$, $(Y, Z)\in (\cS^p, \cM^p)$ for any $p\in (0,1)$.
\end{prop}

\begin{proof}
 It follows from $Y\in \C$ that
 \begin{equation}\label{eq: sup Yp}
  \sup_{0\leq t\leq T} Y_t^p \leq C \pare{1+ \sum_{i=1}^d \sup_{0\leq t\leq T} (X^i_t)^p +\sup_{0\leq t\leq T} \expec\bra{\left. \int_0^T H(s, \uX_s) \,ds \right| \F_t}^p}, \quad \text{ for some constant } C,
 \end{equation}
 where the inequality $(a+b)^p \leq a^p + b^p$ for any $a, b\geq 0$ and $p\in (0,1)$ is used. Recall from Lemma 6.1 in \cite{Briand-et-al} that any 1-dimensional martingale $\{M_t; t\in[0,T]\}$ satisfies
 \[
  \expec[\sup_{0\leq t\leq T} |M_t|^p] \leq \frac{1}{1-p} \expec[|M_T|]^p, \quad \text{ for any } p\in (0,1).
 \]
 Then the previous inequality, Lemma~\ref{lem: local mart Lp}, and \eqref{eq: sup Yp} combined implies that
 \[
  \expec[\sup_{0\leq t\leq T} Y_t^p] \leq \frac{C}{1-p}\pare{1+\sum_{i=1}^d (X_0^i)^p + \expec\bra{\int_0^T H (s, \uX_s) ds}^p} <+\infty, \quad \text{ for any } p\in (0,1).
 \]
 Now recall Lemma 3.1 in \cite{Briand-et-al}. We have from the previous inequality that $Z\in \cM^p$ for any $p\in (0,1)$.
\end{proof}

\section{Viscosity solutions to \eqref{eq: pde}}\label{sec: pde solution}

\subsection{A parabolic boundary value problem}
To show that $u$ and $\ou$, defined in \eqref{eq: u and ou}, are viscosity solutions to \eqref{eq: pde}, we need some preparation first. Given $t\in [0,T]$, a domain $\mathbb{B}_r$ for some $r>0$, and a continuous function $h: \overline{\mathbb{B}_r} \rightarrow \Real$, we consider the BSDE
\begin{equation}\label{eq: Y0x bsde}
 \cY_s = h(X^x_{\tau^{t,x}}) + \int_s^{T-t} \indic_{\{u\leq \tau^{t,x}\}} f(u+t, X^x_u, \cY_u, \cZ_u) \,du -\int_s^{T-t} \cZ_u dB_u,
\end{equation}
where $\tau^{t,x}= \inf\{s\geq 0 \such X^x_s \notin \overline{\mathbb{B}_r}\}\wedge (T-t)$. Since $h(X^x_{\tau^{t,x}})$ is bounded, the previous BSDE admits a unique solution $(\cY^{t,x}, \cZ^{t,x})\in (\cS^\infty, \cM^2)$. Define $(\bY^{t,x}_s, \bZ^{t,x}_s):= (\cY^{t,x}_{(s-t)^+}, \indic_{\{s\geq t\}} \cZ^{t,x}_{s-t})$. They are the unique solution of
\begin{equation}\label{eq: Ytx bsde}
 \bY_s = h(X^{t,x}_{\sigma^{t,x}}) + \int_s^T \indic_{\{t\leq u\leq \sigma^{t,x}\}} f(u, X^{t,x}_u, \bY_u, \bZ_u)\, du -\int_s^T \bZ_u dB_u^t,
\end{equation}
where $\sigma^{t,x}= \inf\{s\geq t\such X^{t,x}_s \notin \overline{\mathbb{B}_r}\} \wedge T$. Now set
\begin{equation}\label{eq: def w}
 w(t,x) := \bY^{t,x}_t = \cY^{t,x}_0, \quad \text{ for } (t,x)\in [0,T]\times(0,\infty)^d.
\end{equation}
Since $\mathbb{S}^+_r$ is regular, $\sigma^{t,x} = t$ for $x\in \mathbb{S}^+_r$, hence $w(t,x) = h(x)$ when $x\in (0,\infty)^d \setminus \mathbb{B}_r^+$.
We claim that
\begin{equation}\label{eq: Y^tx_s}
 w(s, X^{t,x}_s) = \bY^{s, X^{t,x}_s}_s =\bY^{t,x}_s =\cY^{t,x}_{s-t}, \quad \text{ for } s\in [t, \sigma^{t,x}].
\end{equation}
Only the second identity needs a proof. Observe that $X^{s, X^{t,x}_s}_s = X^{t,x}_s$. It follows from the pathwise uniqueness for \eqref{eq: sde X tx} that $\prob(X^{s, X^{t,x}_s}_u = X^{t,x}_u \text{ for any } u\in [s,T])=1$. Then this yields $\sigma^{s, X^{t, x}_s} = \sigma^{t, x}$ for $s\leq \sigma^{t,x}$. The second identity then follows from the uniqueness of solutions to \eqref{eq: Ytx bsde}.

In what follows, we will prove that $w$ is continuous viscosity solution of the following boundary value problem:
\begin{equation}\label{eq: bv}
\begin{array}{ll}
 -\partial_t w - \cL w - f(t, x, w, \nabla w \sigma) =0, \quad & (t,x) \in (0,T) \times \mathbb{B}_r^+,\\
 w(t,x) = h(x), \quad & (t,x)\in (0,T] \times \mathbb{S}_r^+ \cup {T} \times \mathbb{B}_r^+.
\end{array}
\end{equation}
Here no boundary condition is needed on $\mathcal{O}$ because \eqref{H: X non explode} implies that $X^{t,x}$ never reaches $\mathcal{O}$ before $T$. Let us  define what we mean by a continuous viscosity solution of \eqref{eq: bv}.

\begin{defn}
  A continuous function $w: [0,T] \times \overline{\mathbb{B}_r^+} \rightarrow \Real$ is called a \emph{continuous viscosity subsolution} (resp. \emph{supersolution}), if for any $(t,x) \in (0,T) \times \overline{\mathbb{B}_r^+}$, $\phi \in C^{1,2}((0,T) \times (0,\infty)^d)$ such that $(t,x)$ is the local maximum (resp. minimum ) of $w-\phi$, then
 \[
 \begin{array}{ll}
  -\partial_t \phi - \cL \phi - f(t, x, w, \nabla \phi \,\sigma) \leq 0 \, (\text{resp.} \geq 0), & \text{ if } (t,x)\in (0,T)\times \mathbb{B}^+_r,\\
  \min\left\{-\partial_t \phi - \cL \phi - f(t, x, w, \nabla \phi \,\sigma), w(t,x) - h(x)\right\} \leq 0, & \text{ if } (t, x)\in (0,T) \times \mathbb{S}_r^+,\\
  (\text{resp. } \max\left\{\partial_t \phi - \cL \phi - f(t, x, w, \nabla \phi \,\sigma), w(t,x) - h(x)\right\}\geq 0) & \\
  w(T, x) \leq h(x) \, (\text{resp. } w(T, x)\geq h(x)).
 \end{array}
 \]
 A continuous function $w$ is said to be a \emph{continuous viscosity function} if it is both viscosity sub- and supersolution.
\end{defn}

Since points on $\mathbb{S}_r^+$ are regular, the following result can be viewed as the parabolic analogue of Proposition 6.3 and Theorem 6.5 in \cite{Darling-Pardoux}, where a similar result has been proved for an elliptic boundary value problem. 

\begin{prop}\label{prop: w cont vis}
 Suppose that \eqref{H: f lip z} - \eqref{H: f nonneg}, \eqref{H: X non explode}, and \eqref{H: nondeg} hold. Then $w$ is a continuous viscosity solution to \eqref{eq: bv}.
\end{prop}

\subsection{Existence of viscosity solutions of \eqref{eq: pde}}\label{subsec: existence vis}
Now choosing different $h$ in \eqref{eq: Y0x bsde} and \eqref{eq: Ytx bsde}, we can construct approximating sequences for $u$ and $\overline{u}$. For each $n$, choose $r=n$, we rename $w$ in \eqref{eq: def w} as $u_n$ when $h=g$, and $\ou_n$ when $h= g_n$. Both $u_n$ and $\overline{u}_n$ are defined on $[0,T]\times (0,\infty)^d$. Solutions to \eqref{eq: Y0x bsde} and \eqref{eq: Ytx bsde} are denoted as $(\cY^{n, t, x}, \cZ^{n, t, x})$ and $(\bY^{n, t, x}, \bZ^{n, t, x})$ respectively when $h= g$; and $(\overline{\cY}^{n, t, x}, \overline{\cZ}^{n, t, x})$ and $(\overline{\bY}^{n, t, x}, \overline{\bZ}^{n, t, x})$ respectively when $h= g_n$. Then $u_n(t, x) = Y^{n, t, x}_t =\cY^{n, t, x}_0$, for $(t,x) \in [0,T]\times (0, \infty)^d$, and a similar identity holds for $\ou_n$ as well. Note that $\mathbb{B}^+_r$ increases to $(0,\infty)^d$ as $n\rightarrow \infty$, it follows from the definition of $u$ and $\ou$ in \eqref{eq: u and ou} and the construction of $\cY^{t,x}$ and $\overline{\cY}^{t,x}$ in \eqref{eq: bsde cY} and Theorem \ref{thm: bsde existence} that
\[
 \uparrow \lim_{n\rightarrow \infty} u_n(t,x) = \uparrow \lim_{n\rightarrow \infty} \cY^{n, t, x}_0 = \cY^{t, x}_0 = u(t,x), \quad \text{ for } (t,x)\in [0,T] \times (0,\infty)^d,
\]
where the second identity holds for $t=T$ thanks to $\uparrow \lim_{n\rightarrow \infty} g_n(x) = g(x)$.
A similar statement holds for $\ou_n$ and $\ou$ as well. On the other hand, Proposition \ref{prop: w cont vis} implies that $u_n$ (resp. $\ou_n$) is a continuous viscosity solution to the boundary value problem \eqref{eq: bv} when the boundary condition is $g$ (resp. $g_n$).

Before using $\{u_n\}_{n\geq 0}$ and $\{\ou_n\}_{n\geq 0}$ to prove that both $u$ and $\ou$ solves \eqref{eq: pde} in the viscosity sense defined in Definition \ref{def: solution}, we recall \emph{half-relaxed upper} and \emph{lower limits} of $\{u_n\}_{n\geq 0}$:
\[
\begin{split}
u^{U}(t,x)&:= \limsup_{n\rightarrow \infty} \left\{u_m(t', x') \such m\geq n, (t', x') \in (0,T) \times (0,\infty)^d, \text{ and } |t'-t|+|x'-x| \leq 1/n\right\},\\
u^{L}(t,x)&:= \liminf_{n\rightarrow \infty} \left\{u_m(t', x') \such m\geq n, (t', x') \in (0,T) \times (0,\infty)^d, \text{ and } |t'-t|+|x'-x| \leq 1/n\right\}.
\end{split}
\]
The half-relaxed upper and lower limits $\ou^U$ and $\ou^L$ are defined analogously for $\{\ou_n\}_{n\geq 0}$.

Since $\{u_n\}_{n\geq 0}$ is a nondecreasing sequence of continuous functions, the following orders among $u, u^U, u^L, u^*$  and $u_*$ hold. The same order also holds for functions associated to $\ou$ as well.

\begin{lem}\label{lem: order four funs}
 For any $(t,x) \in [0,T] \times (0,\infty)^d$,
 \[
  u(t,x) = u^L(t,x) = u_*(t,x) \leq u^*(t,x) = u^U(t,x).
 \]
\end{lem}

\begin{proof}
This relationship has been applied in \cite{Popier}. But no reference or proof is given there. For the reader's convenience, we present a short proof here.

\underline{$u= u_*$}: Since $u$ is the supremum of continuous functions $\{u_n\}_{n\geq 0}$, $u$ is lower-semicontinuous. Recall that $u_*$ is the largest lower-semicontinuous function dominated by $u$. Hence $u=u_*$.

\underline{$u^L \leq u$}: Since $u$ is lower-semicontinuous, there exists a sequence $(t_n, x_n) \in B_{(t,x)}(1/n)$ such that $\lim_{n\rightarrow \infty} u(t_n, x_n) = u(t,x)$. Here $B_{(t,x)}(r):= \{(t',x')\in (0,T)\times (0,\infty)^d \such |t'-t|+|x'-x| \leq r\}$. Since $\{u_n\}_{n\geq 0}$ is nondecreasing,
\[
 \inf_{(t', x') \in B_{(t,x)}(1/n)\,;\, m\geq n} u_m(t', x') \leq u_n(t_n, x_n) \leq u(t_n, x_n).
\]
The claim then follows from sending $n\rightarrow \infty$ in the previous inequalities.

\underline{$u^L \geq u$}: For any $n\geq N \geq 0$,
\[
 u^L(t, x) \geq \inf_{(t', x')\in B_{(t,x)}(1/n)\,;\, m\geq n} u_m(t', x') \geq \inf_{(t', x') \in B_{(t,x)}(1/n)} u_N(t', x'),
\]
where the second inequality holds since $\{u_n\}_{n\geq 0}$ is nondecreasing. Now, sending $n\rightarrow \infty$ and using the continuity of $u_N$, we obtain from the previous inequalities that $u^L(t,x) \geq u_N(t, x)$. The claim then follows after sending $N\rightarrow \infty$.

\underline{$u^U \leq u^*$}: Let $(t_n, x_n)$ be a sequence converging to $(t,x)$ such that $\lim_{n\rightarrow \infty} u_n(t_n, x_n) = u^U(t, x)$. Since $\{u_n\}_{n\geq 0}$ is a nondecreasing sequence,
$
 u_n(t_n, x_n) \leq u(t_n, x_n).
$
Sending $n\rightarrow \infty$, the claim follows from the upper semicontinuity of $u^*$.

\underline{$u^*\leq u^U$}: For any $\epsilon>0$, these exists a sufficiently large $N$, such that
\[
 \epsilon + u^U(t, x) \geq \sup_{(t', x') \in B_{(t,x)}(1/N)\,;\, n \geq N} u_n(t', x') \geq u_n(\wt{t}, \wt{x}), \quad \text{ for any } n\geq N \text{ and } (\wt{t}, \wt{x}) \in B_{(t,x)}(1/n).
\]
Since $\{u_n\}_{n\geq 0}$ is nondecreasing, the previous inequality yields $\epsilon + u^U(t, x) \geq u(\wt{t}, \wt{x})$. Now the claim follows from first sending $n\rightarrow \infty$ then $\epsilon \rightarrow 0$.
\end{proof}

From the previous lemma and the definition of $u$,  we have $u^*(T,x) \geq u_*(T,x) = u(T, x) =g(x)$ for any $x\in (0,\infty)^d$. In what follows, we will prove the converse inequality. The same statement holds for $\overline{u}^*$ as well.

\begin{lem}\label{lem: u* terminal}
 Let \eqref{H: g linear growth} - \eqref{H: nondeg} hold. Then
 $u^*(T,x) \leq g(x)$ for  $x\in (0,\infty)^d$.
\end{lem}

\begin{proof}
 It suffices to prove the statement for $u^U$, since $u^U = u^*$. Take any sequence $\{(t_n, x_n)\}_{n\geq 0}$ converging to $(T,x)$. Without loss of generality, we can assume all $x_n \in D$ for a bounded domain $D\subset (0,\infty)^d$ containing $x$. Recall that $u_m(t_n, x_n) =\cY^{m, t_n, x_n}_0$, where $(\cY^{m}, \cZ^{m}) \in (\cS^\infty, \cM^2)$ (the superscript $(t_n, x_n)$ is omitted for simplicity of notation) is the unique solution of the following BSDE:
 \[
  \cY^m_s = g_m(X^{x_n}_{\tau^{t_n, x_n}_m}) + \int_s^{T-t_n} \indic_{\{u \leq \tau^{t_n, x_n}_m\}} f(u+t_n, X^{x_n}_u, \cY^m_u, \cZ^m_u) \,du - \int_s^{T-t_n} \cZ^m_u \,dB_u,
 \]
 where $\tau^{t_n, x_n}_m = \inf\{u\geq 0 \such X^{x_n}_u \notin \overline{\mathbb{B}_m}\}\wedge(T-t_n)$.
 Choosing $s=0$ and taking expectation in the last equation, we obtain
 \begin{equation}\label{eq: Ym0 est}
  \cY^m_0 = \expec\bra{g_m(X^{x_n}_{\tau^{t_n, x_n}_m})} + \expec\bra{\int_0^{\tau^{t_n, x_n}_m}f(u+ t_n, X^{x_n}_u, \cY^m_u, \cZ^m_u) \,du}.
 \end{equation}

 Let us estimate individual terms on the right hand side of the previous identity.
 Assumptions (H3), (H6) and Lemma \ref{lem: Y bound} combined implies that
 \[
  0\leq f(u+t_n, X^{x_n}_u, \cY^m_u, \cZ^m_u) \leq \wt{K}(1+ \uX^{x_n}_u) + \mu \cY^m_u \leq C(1+\uX^{x_n}_u),
 \]
 where $C$ is a constant depending on $K, \wt{K}, \mu, T$, but not $n$ and $m$. As a result,
 \[
  \lim_{n\rightarrow \infty} \expec\bra{\int_0^{\tau^{t_n, x_n}_m}f(u+ t_n, X^{x_n}_u, \cY^m_u, \cZ^m_u) \,du} \leq \lim_{n\rightarrow \infty} C\int_0^{T-t_n} 1+ \expec[\uX^{x_n}_u] \,du \leq \lim_{n\rightarrow \infty} C(T-t_n) (1+\underline{x_n}) =0,
 \]
 where the second inequality follows from the supermartingale property of $\uX^{x_n}$. On the other hand,
 \begin{equation}\label{eq: gm exp upper bound}
 \begin{split}
  \expec\bra{g_m(X^{x_n}_{\tau^{t_n, x_n}_m})} \leq \expec\bra{g(X^{x_n}_{\tau^{t_n, x_n}_m})} &= K \pare{1+ \expec[\uX^{x_n}_{\tau^{t_n, x_n}_m}]} - \expec\bra{\overline{g}(X^{x_n}_{\tau^{t_n, x_n}_m})} \\
  & = K(1+\underline{x_n}) - \expec\bra{\overline{g}(X^{x_n}_{\tau^{t_n, x_n}_m})}.
 \end{split}
 \end{equation}
 Recall that coefficients in \eqref{eq: sde X} is locally Lipschitz. It then follows from the continuity of stochastic flow (see Theorem 5.38 in \cite{Protter}) that $\prob-\lim_{n\rightarrow \infty} \sup_{s\in [0,T]} |X^{x_n}_s - X^x_s| =0$. As a result,
 \[
  |X^{x_n}_{\tau^{t_n, x_n}_m} - x| \leq \sup_{s\in[0,T]}|X^{x_n}_s - X^x_s| + |X^{x}_{\tau^{t_n, x_n}_m} -x| \rightarrow 0 \quad \text{ under } \prob,
 \]
 as $n\rightarrow \infty$. Going back to \eqref{eq: gm exp upper bound}, we have from Fatou's lemma that
 \[
  \limsup_{m \geq n, n\rightarrow \infty} \expec\bra{g_m(X^{x_n}_{\tau^{t_n, x_n}_m})} \leq K(1+\underline{x}) - \liminf_{m\geq n, n\rightarrow \infty} \expec\bra{\overline{g}(X^{x_n}_{\tau^{t_n, x_n}_m})} \leq K(1+\underline{x}) - \og(x) = g(x).
 \]
 Therefore, we conclude from \eqref{eq: Ym0 est} that
 $\limsup_{m\geq n\,;\, n\rightarrow \infty} u_m(t_n, x_n) \leq g(x)$. The statement then follows since the choice of $(t_n, x_n)_{n\geq 0}$ is arbitrary.
\end{proof}

Let $D\subset \Real^d$ be locally compact and $D_T = (0,T) \times D$. We recall \emph{parabolic semijets} $\mathcal{P}^{2, \pm}$ from \cite{user-guide}. The proof of Theorem \ref{thm: pde existence} needs the following stability property of parabolic semijets. This result is a straight forward extension of Proposition 4.2 in \cite{user-guide} to its parabolic analogue.

\begin{lem}\label{lem: vis stability}
 Let $v$ be a upper semi-continuous function on $D_T$, $(t,z)\in D_T$, and $(a, p, W) \in \mathcal{P}^{2, +} v(t,z)$. Suppose also that $v_n$ is a sequence of upper semi-continuous functions on $D_T$ such that
 \begin{enumerate}
 \item[(i)] there exists $(t_n, x_n) \in D_T$ such that $(t_n, x_n, v_n(t_n, x_n)) \rightarrow (t, z, v(t,z))$,
 \item[(ii)] if $(s_n, z_n)\in D_T$ and $(s_n, z_n) \rightarrow (s, x)\in D_T$, then $\limsup_{n\rightarrow \infty} v_n(s_n, z_n) \leq v(s, x)$.
 \end{enumerate}
 Then there exists $(\hat{t}_n, \hat{x}_n)\in D_T$, $(a_n, p_n, W_n) \in \mathcal{P}^{2, +} v_n(\hat{t}_n, \hat{x}_n)$ such that
 \[(\hat{t}_n, \hat{x}_n, v_n(\hat{t}_n, \hat{x}_n), a_n, p_n, W_n) \rightarrow (t, z, v(t,z), a, p, W).\]
\end{lem}

Now we are ready to prove that both $u$ and $\ou$ are viscosity solutions to \eqref{eq: pde}.
\begin{proof}[Proof of Theorem \ref{thm: pde existence}]
 We have already seen $u(t,x) > \ou(t,x)$ for $(t,x)\in [0,T)\times (0, \infty)^d$. It only remains to show that both $u$ and $\ou$ are viscosity solutions of at most linear growth. We will only prove the statement for $u$. The statement for $\ou$ can be proved similarly.
 First, Theorem \ref{thm: pde existence} and \eqref{H: linear H} combined implies that $u(t, x) = \cY^{t,x}_0 \leq C(1+\underline{x})$ where $C$ is a constant depending on $\mu, K, \wt{K}$, and $T$. Second, $u^*(T, x)\leq g(x)$ has already been proved in Lemma \ref{lem: u* terminal}. Lastly, for any $(t,x) \in (0,T) \times (0,\infty)^d$ and $(a, p, W) \in \mathcal{P}^{2, +}u^*(t,x)$, since $u^* = u^U$, we want to show
 \[
  -a - \frac12 Tr(\sigma \sigma'(x) W) - f(t, x, u^U(t, x), p) \leq 0.
 \]
 Since there is a sufficiently large $\mathbb{B}^+_n$ such that $(t, x) \in (0,T)\times \mathbb{B}^+_n$, the previous inequality follows directly from Proposition \ref{prop: w cont vis} and Lemma \ref{lem: vis stability}. Similar argument shows that $u$ is also a supersolution.
\end{proof}

\subsection{Uniqueness of viscosity solutions of \eqref{eq: pde}}\label{subsec: unique vis}

To prove the comparison result, let us first present the following lemma, which is similar to Lemma 3.7 in \cite{Barles-et-al}.
\begin{lem}\label{lem: u-v vis sub}
 Let $u$ be a subsolution, $v$ be a supersolution of \eqref{eq: pde},  and both $u$ and $v$ be locally bounded in $[0,T] \times (0, \infty)^d$. Then $w:= u-v$ is a viscosity subsolution of
 \begin{equation}\label{eq: w pde}
  -\partial_t w - \cL w - \mu |w| - b\, |\nabla w \,\sigma| =0, \quad \text{ in } [0,T) \times (0,\infty)^d.
 \end{equation}
\end{lem}

\begin{proof}
 The proof is essentially the same with the proof of Lemma 3.7 in \cite{Barles-et-al}, except several points which we are going to emphasize as follows. Let us follow the notation in \cite{Barles-et-al}.

 First, fix $(t_0, x_0)\in (0,T) \times (0,\infty)^d$ and $\phi(t, x)\in C^{1,2}((0,T)\times (0,\infty)^d)$ such that $w^*-\phi$ attains its strict global maximum at $(t_0, x_0)$. We double the variables and introduce an auxiliary function
 \[
  \psi_{\epsilon, \alpha}(t, x, s, y) = u^*(t, x) - v_*(s, y) - \frac{|x-y|^2}{\epsilon^2} - \frac{|t-s|^2}{\alpha^2} - \phi(t,x),
 \]
 where $\epsilon, \alpha$ are positive parameters which tend to zero. Fix a compact domain $B \subset (0,\infty)^d$ which contains $x_0$. Since $\psi_{\epsilon, \alpha}$ is upper semi-continuous and is bounded from above on $([0,T] \times B)^2$, the maximum of $\psi_{\epsilon, \alpha}$ on $([0, T] \times B)^2$ is attained at a point $(\overline{t}, \overline{x}, \overline{s}, \overline{y}) \in ([0,T]\times B)^2$. We have dropped the dependence of $\overline{t}, \overline{x}, \overline{s}$, and $\overline{y}$ in $\epsilon$ and $\alpha$ for simplicity of notation. We claim that
 \begin{enumerate}
  \item[(i)] $(\overline{t}, \overline{x}), (\overline{s}, \overline{y}) \rightarrow (t_0, x_0)$ as $\epsilon, \alpha \rightarrow 0$,
  \item[(ii)] $|\overline{x} - \overline{y}|^2/\epsilon^2$ and $|\overline{t} - \overline{s}|^2/\alpha^2$ are bounded and tend to zero as $\epsilon, \alpha \rightarrow 0$.
 \end{enumerate}
 Indeed, since $\psi_{\epsilon, \alpha}$ attains its maximum at $(\overline{t}, \overline{x}, \overline{s}, \overline{y})$,
 \begin{equation}\label{eq: max w-phi}
  u^*(\overline{t}, \overline{x}) - v_*(\overline{s}, \overline{y}) - \frac{|\overline{x} -\overline{y}|^2}{\epsilon^2} - \frac{|\overline{t}-\overline{s}|^2}{\alpha^2} -\phi(\overline{t}, \overline{x}) \geq u^*(t_0, x_0) - v_*(t_0, x_0) - \phi(t_0, x_0).
 \end{equation}
 The previous inequality yields
 \[
  |\overline{x} - \overline{y}|^2 + |\overline{t} -\overline{s}|^2 \leq 2\, (\epsilon^2 \vee \alpha^2) \max_{(t,x,s,y) \in ([0,T] \times B)^2} |u^*(t, x)-v_*(s, y) - \phi(t,x)|.
 \]
 Send $\epsilon, \alpha \rightarrow 0$ in the previous inequality. It follows that $\overline{x}$ and $\overline{y}$ converge to the same point, say $\wt{x}$, meanwhile $\overline{t}$ and $\overline{s}$ converge to $\wt{t}$. Then sending $\epsilon, \alpha \rightarrow 0$ on the left side of \eqref{eq: max w-phi} and using the upper semi-continuity of $u^*-v_*-\phi$, we obtain
 \[
  u^*(\wt{t}, \wt{x}) - v_*(\wt{t}, \wt{x}) - \phi(\wt{t}, \wt{x}) - \lim_{\epsilon, \alpha \rightarrow 0} \pare{\frac{|\overline{x} -\overline{y}|^2}{\epsilon^2} + \frac{|\overline{t}-\overline{s}|^2}{\alpha^2}} \geq u^*(t_0, x_0) - v_*(t_0, x_0) - \phi(t_0, x_0).
 \]
 Since $u^*-v_*-\phi$ attains its strict global maximum at $(t_0, x_0)$, both claims follow from the previous inequality. Now apply Theorem 8.3 in \cite{user-guide} to obtain two triplets $(\overline{a} + \partial_t \phi(\overline{t}, \overline{x}), \overline{p} + \nabla \phi(\overline{t}, \overline{x}), M) \in \overline{\mathcal{P}}^{2, +} u^*(\overline{t}, \overline{x})$ and $(\overline{a}, \overline{p}, N) \in \overline{\mathcal{P}}^{2, -} v_*(\overline{s}, \overline{y})$ respectively, and write down two inequalities that these triplets satisfy. When we estimate the difference between these two inequalities, since $\overline{x}, \overline{y} \in B$, we can use the local Lipschitz continuity of $\sigma$ on $B$:
 \[
  Tr(\sigma \sigma' (\overline{x}) M) - Tr(\sigma \sigma' (\overline{y}) N) \leq C_B \frac{|\overline{x} - \overline{y}|^2}{\epsilon^2} + Tr(\sigma \sigma'(\overline{x}) \phi(\overline{t}, \overline{x})),
 \]
 where $C_B$ depends on the Lipschitz constant of $\sigma$ on $B$.

 Second, since \eqref{H: f lip z str} is assumed,
 \[
  f(\cdot, \overline{x}, \cdot, (\overline{p} + \nabla \phi(\overline{t}, \overline{x})) \sigma(\overline{x})) - f(\cdot, \overline{x}, \cdot, \overline{p} \sigma(\overline{y})) \leq b(\overline{x}) \left|\overline{p} \,(\sigma(\overline{x}) - \sigma(\overline{y})) + \nabla \phi(\overline{t}, \overline{x}) \sigma(\overline{x})\right|,
 \]
 where $\overline{p} = 2(\overline{x} - \overline{y}) / \epsilon^2$. Using the Lipschitz continuity of $\sigma$ on $B$, the right side of the previous inequality converges to $b(x_0) \left|\nabla \phi(t_0, x_0) \sigma(x_0)\right|$ as $\epsilon, \alpha\rightarrow 0$. The rest argument is the same with that in Lemma 3.7 in \cite{Barles-et-al}.
\end{proof}

In what follows we are going to construct a supersolution to \eqref{eq: w pde}, using the Lyapunov function $\Psi$ in (H9).
\begin{lem}\label{lem: w pde supersoln}
 Assume that the strict positive function $\Psi$ in (H9) exists. Then $\Phi(t, x):= e^{L(T-t)} \Psi(x)$, for sufficiently large $L$, satisfies
 \[
  -\partial_t \Phi - \cL \Phi - \mu \Phi - b |\nabla \Phi \sigma| >0, \quad \text{ in } [0,T]\times (0,\infty)^d.
 \]
\end{lem}
\begin{proof}
 Since $\Psi(x) \rightarrow \infty$ as $x\rightarrow \overline{x} \in \mathcal{O}$ or $|x| \rightarrow \infty$, $\Psi$ being strictly positive implies that $m:= \min_{x\in (0,\infty)^d} \Psi(x) >0$. It then follows from \eqref{H: Lpsi} and \eqref{H: Psi der bound} that
 \[
 \begin{split}
  & - \partial_t \Phi(t,x) - \cL \Phi(t,x) - \mu \Phi(t,x) - b(x) \left|\nabla_x \Phi(t, x) \sigma(x)\right|\\
  & \hspace{1cm} \geq L \Phi(t, x) - \lambda e^{L(T-t)} - \lambda \Phi(t,x) - \mu \Phi(t,x) - c \Phi(t,x).
 \end{split}
 \]
 Since $\Phi(t,x) \geq m$ on $[0,T] \times (0,\infty)^d$, one can choose sufficiently large $L$ such that the right side of the previous inequality is strictly positive on $[0,T] \times (0,\infty)^d$.
\end{proof}

\begin{proof}[Proof of Theorem \ref{thm: pde uniqueness}]
 We are going to show that $w= u-v$ satisfies $w(t,x) \leq \alpha \Phi(t,x)$ in $[0,T] \times (0,\infty)^d$ for any $\alpha>0$. Sending $\alpha$ to zero, we obtain $u \leq v$ on $[0,T]\times (0,\infty)^d$.

 To prove the claim, let us consider $M(t,x) := \pare{w^*(t,x) - \alpha \Phi(t,x)} e^{\mu t}$. Since $u$ and $v$ are nonnegative and bounded from above by $C(1+\underline{x})$ for some constant $C$, moreover \eqref{H: psi close O} and \eqref{H: Psi close inf} imply that $\Phi(t,x) \rightarrow \infty$ as $x\rightarrow \overline{x} \in \mathcal{O}$ and $\lim_{|x|\rightarrow \infty} \Phi(x) / \underline{x} =\infty$, then there exists a compact domain $B\subset (0,\infty)^d$ such that $M(t,x) <0$ for $(t,x) \in [0,T] \times B^c$. On the other hand, since $M$ is upper semi-continuous, it attains its maximum in $[0,T]\times B$ at a point, say $(t_0, x_0)$. We can assume that $M(t_0, x_0) > 0$, otherwise $M \leq 0$ on $[0,T]\times (0,\infty)^d$ and we are done. As a result $(t_0, x_0)$ is the global maximum point of $M$ on $[0,T]\times (0,\infty)^d$ and $w^*(t_0, x_0)>0$, which implies $t_0< T$ thanks to $u^*(T, x) \leq g(x) \leq v_*(T, x)$.

 From the maximum point property, we obtain that
 \[
  w^*(t,x) - \alpha \Phi(t,x) \leq \pare{w^*(t_0, x_0) - \alpha \Phi(t_0, x_0)} e^{\mu (t_0 - t)}, \quad \text{ for any } (t,x)\in [0,T]\times (0,\infty)^d.
 \]
 This inequality implies that $w^* - \phi$ attains its global maximum point at $(t_0, x_0)$, where
 \[
  \phi(t, x) = \alpha \Phi(t,x) + \pare{w^*(t_0, x_0) - \alpha \Phi(t_0, x_0)} e^{\mu (t_0 - t)}.
 \]
 Since $w$ is a subsolution to \eqref{eq: w pde}, we have
 \[
  -\partial_t \phi(t_0, x_0) - \cL \phi(t_0, x_0) - \mu w^*(t_0, x_0) - b(x_0) \left|\nabla \phi(t_0, x_0) \sigma(x_0)\right| \leq 0.
 \]
 But the left side of this inequality is
 \[
  \alpha\bra{-\partial_t \Phi(t_0, x_0) -\cL \Phi(t_0, x_0) - \mu \Phi(t_0, x_0) - b(x_0) \left|\nabla \Phi(t_0, x_0) \sigma(x_0)\right|}.
 \]
 We then obtain a contradiction with Lemma \ref{lem: w pde supersoln}.
\end{proof}

\bibliographystyle{siam}
\bibliography{biblio}

\end{document}